\documentclass[a4paper,12pt]{article}
\usepackage[english]{babel}
\usepackage{graphicx}
\usepackage{amssymb}
\usepackage{amsmath}
\usepackage{geometry}
\usepackage{hyperref}
\hypersetup{
colorlinks=true,
linkcolor=black,
citecolor=black,
urlcolor=black
}

\begin{document}

\renewcommand{\r}{\mathbb R}
\newcommand{\eqd}{\stackrel{d}{=}}

\renewcommand{\refname}{References}

\title{\vspace{-1.8cm}
On the asymptotic approximation to the probability distribution of
extremal precipitation}

\author{
V.\,Yu.~Korolev\textsuperscript{1},
A.\,K.~Gorshenin\textsuperscript{2} }

\date{}

\maketitle

\footnotetext[1]{Faculty of Computational Mathematics and
Cybernetics, Lomonosov Moscow State University, Russia; Institute of
Informatics Problems, Federal Research Center ``Computer Science and
Control'' of Russian Academy of Sciences, Russia; Hangzhou Dianzi
University, China; \url{vkorolev@cs.msu.su}}

\footnotetext[2]{Institute of Informatics Problems, Federal Research
Center ``Computer Science and Control'' of Russian Academy of
Sciences, Russia; \url{agorshenin@frccsc.ru}}

{\bf Abstract.} Based on the negative binomial model for the
duration of wet periods measured in days~\cite{Gulev}, an asymptotic
approximation is proposed for the distribution of the maximum daily
precipitation volume within a wet period. This approximation has the
form of a scale mixture of the Fr{\'e}chet distribution with the
gamma mixing distribution and coincides with the distribution of a
positive power of a random variable having the Snedecor--Fisher
distribution. The proof of this result is based on the
representation of the negative binomial distribution as a mixed
geometric (and hence, mixed Poisson) distribution \cite{Korolev2017}
and limit theorems for extreme order statistics in samples with
random sizes having mixed Poisson distributions
\cite{KorolevSokolov2008}. Some analytic properties of the obtained
limit distribution are described. In particular, it is demonstrated
that under certain conditions the limit distribution is mixed
exponential and hence, is infinitely divisible. It is shown that
under the same conditions the limit distribution can be represented
as a scale mixture of stable or Weibull or Pareto or folded normal
laws. The corresponding product representations for the limit random
variable can be used for its computer simulation. Several methods
are proposed for the estimation of the parameters of the
distribution of the maximum daily precipitation volume. The results
of fitting this distribution to real data are presented illustrating
high adequacy of the proposed model. The obtained mixture
representations for the limit laws and the corresponding asymptotic
approximations provide better insight into the nature of mixed
probability (``Bayesian'') models.

\smallskip


\section{Introduction}

In most papers dealing with the statistical analysis of
meteorological data available to the authors, the suggested
analytical models for the observed statistical regularities in
precipitation are rather ideal and inadequate. For example, it is
traditionally assumed that the duration of a wet period (the number
of subsequent wet days) follows the geometric distribution (for
example, see~\cite{Zolina2013}) although the goodness-of-fit of this
model is far from being admissible. Perhaps, this prejudice is based
on the conventional interpretation of the geometric distribution in
terms of the Bernoulli trials as the distribution of the number of
subsequent wet days (``successes'') till the first dry day
(``failure''). But the framework of Bernoulli trials assumes that
the trials are independent whereas a thorough statistical analysis
of precipitation data registered in different points demonstrates
that the sequence of dry and wet days is not only independent, but
it is also devoid of the Markov property so that the framework of
Bernoulli trials is absolutely inadequate for analyzing
meteorological data.

It turned out that the statistical regularities of the number of
subsequent wet days can be very reliably modeled by the negative
binomial distribution with the shape parameter less than one. For
example, in~\cite{Gulev} we analyzed meteorological data registered
at two geographic points with very different climate: Potsdam
(Brandenburg, Germany) with mild climate influenced by the closeness
to the ocean with warm Gulfstream flow and Elista (Kalmykia, Russia)
with radically continental climate. The initial data of daily
precipitation in Elista and Potsdam are presented on Figures~1a and
1b, respectively. On these figures the horizontal axis is discrete
time measured in days. The vertical axis is the daily precipitation
volume measured in centimeters. In other words, the height of each
``pin'' on these figures is the precipitation volume registered at
the corresponding day (at the corresponding point on the horizontal
axis).

\renewcommand{\figurename}{\rm{Fig.}}

\begin{figure}[h]
\begin{minipage}[h]{0.49\textwidth}
\center{\includegraphics[width=\textwidth,
height=0.6\textwidth]{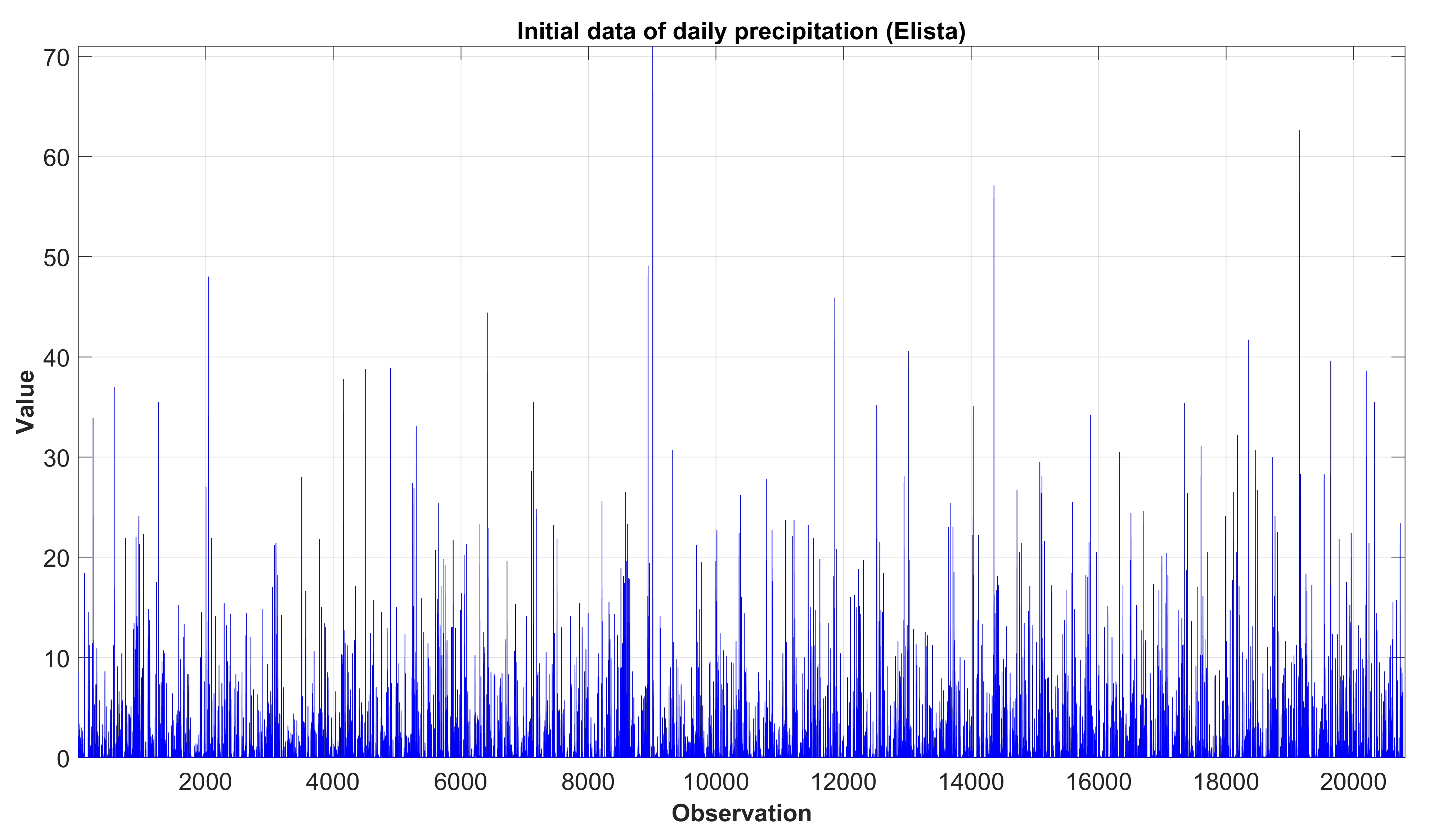}
\\a)}
\end{minipage}
\hfill
\begin{minipage}[h]{0.49\textwidth}
\center{\includegraphics[width=\textwidth]{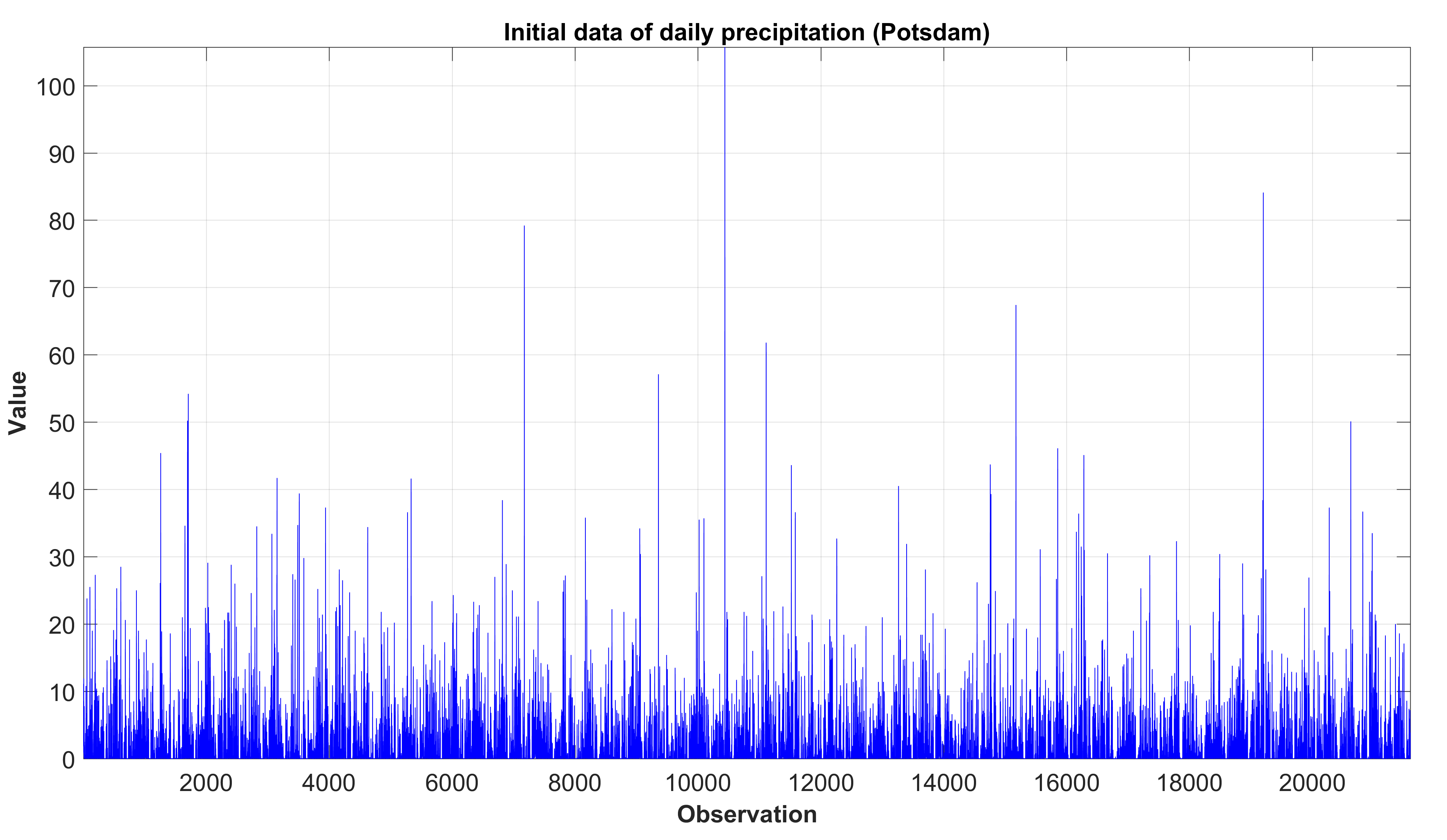} \\
b)}
\end{minipage}
\label{Data} \caption{The initial data of daily precipitation in
Elista (a) and Potsdam (b).}
\end{figure}

In order to analyze the statistical regularities of the duration of
wet periods this data was rearranged as shown on Figures~2a and 2b.

\begin{figure}[h]
\begin{minipage}[h]{0.49\textwidth}
\center{\includegraphics[width=\textwidth,
height=0.6\textwidth]{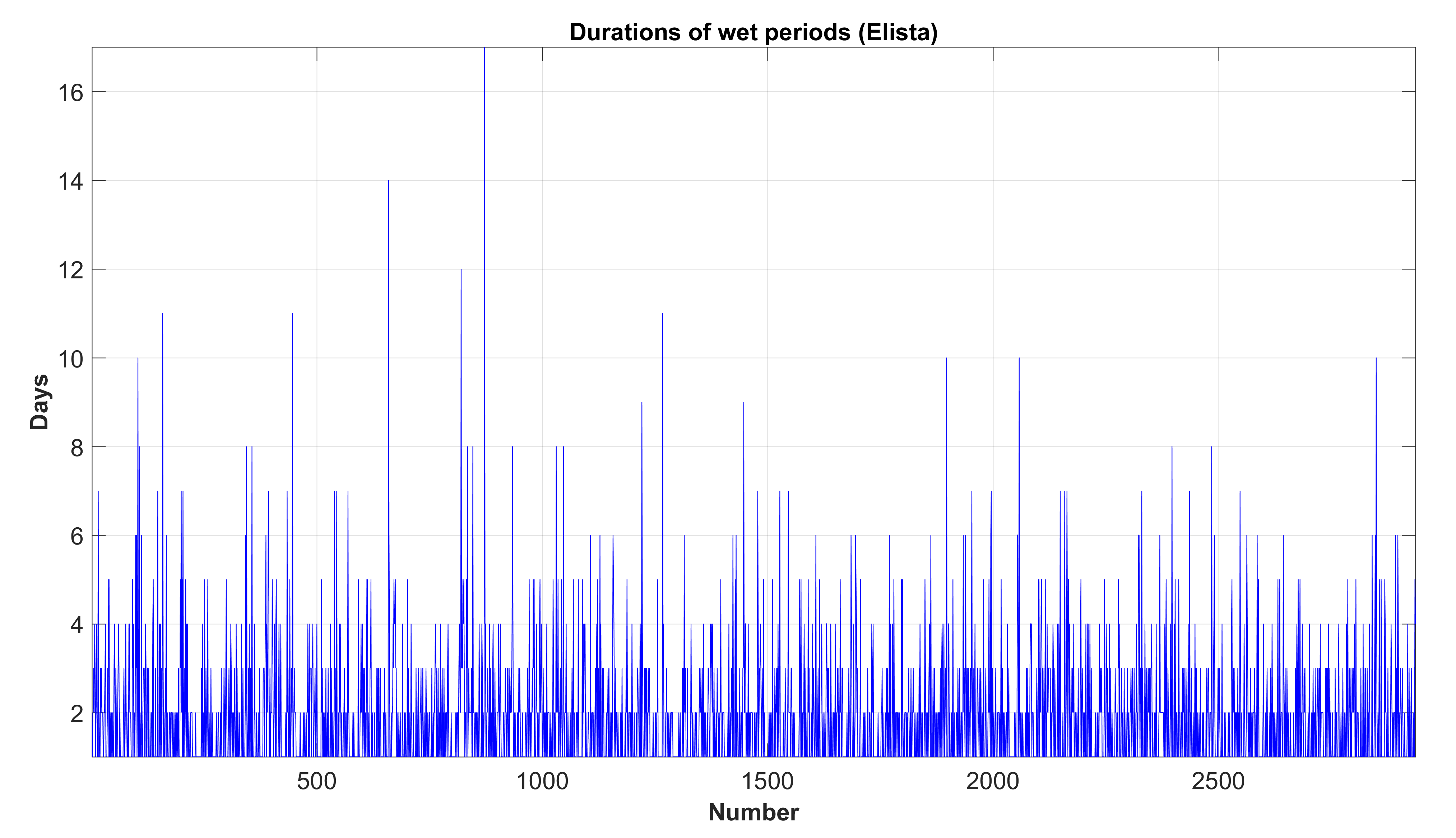}
\\a)}
\end{minipage}
\hfill
\begin{minipage}[h]{0.49\textwidth}
\center{\includegraphics[width=\textwidth]{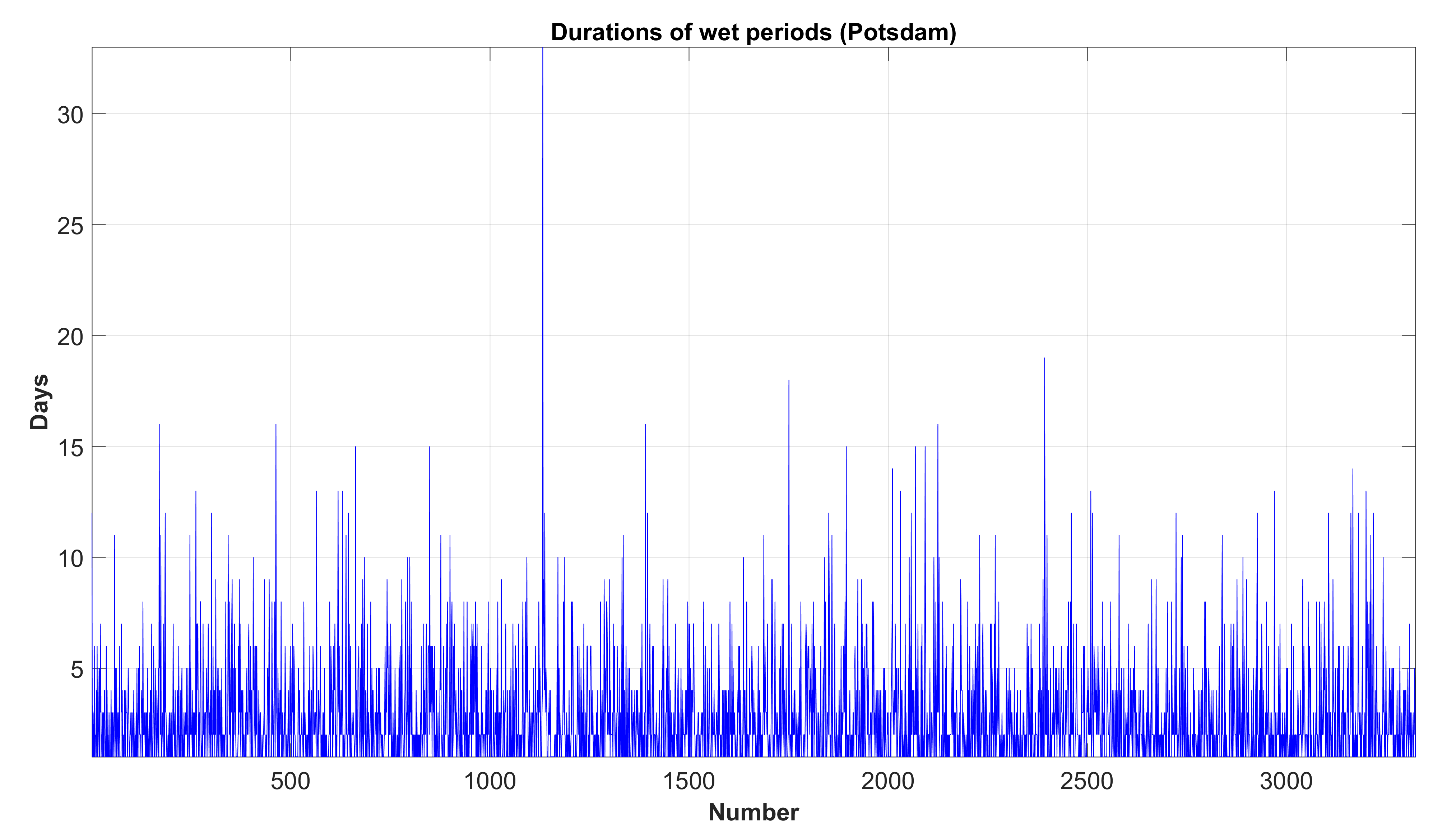} \\
b)}
\end{minipage}
\label{WetPeriod} \caption{The durations of wet periods in Elista
(a) and Potsdam (b).}
\end{figure}

On these figures the horizontal axis is the number of successive wet
periods. It should be mentioned that directly before and after each
wet period there is at least one dry day, that is, successive wet
periods are separated by dry periods. On the vertical axis there lie
the durations of wet periods. In other words, the height of each
``pin'' on these figures is the length of the corresponding wet
period measured in days and the corresponding point on the
horizontal axis is the number of the wet period.

The samples of durations in both Elista and Potsdam were assumed
homogeneous and independent. It was demonstrated that the
fluctuations of the numbers of successive wet days with very high
confidence fit the negative binomial distribution with shape
parameter less than one (also see~\cite{Gorshenin2017}). Figures~3a
and~3b show the histograms constructed from the corresponding
samples of duration periods and the fitted negative binomial
distribution. In both cases the shape parameter $r$ turned out to be
less than one. For Elista $r=0.876$, $p=0.489$, for Potsdam
$r=0.847$, $p=0.322$.

\begin{figure}[h]
\begin{minipage}[h]{0.5\textwidth}
\center{\includegraphics[width=\textwidth,
height=0.6\textwidth]{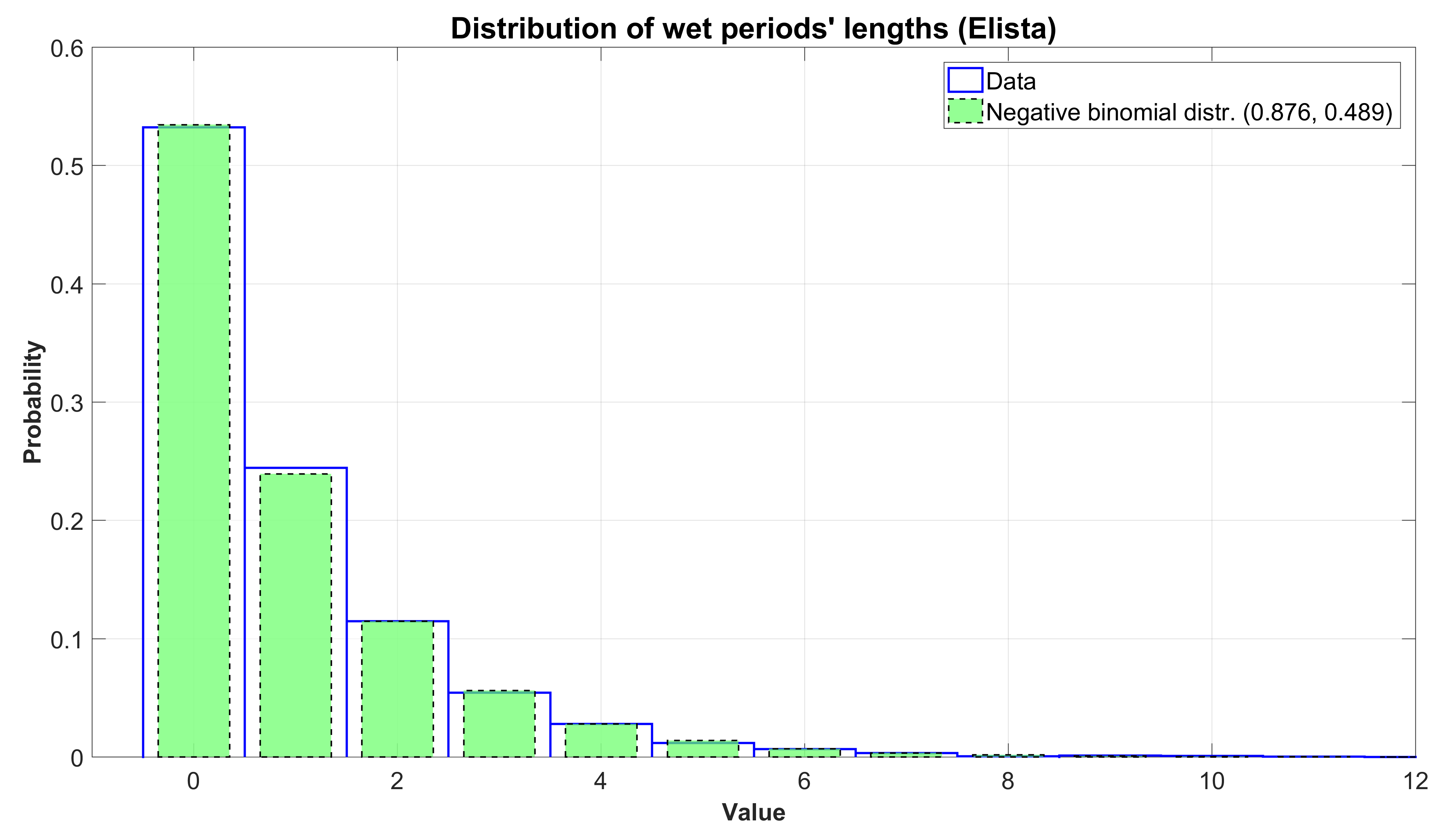}
\\a)}
\end{minipage}
\hfill
\begin{minipage}[h]{0.5\textwidth}
\center{
\includegraphics[width=\textwidth,
height=0.6\textwidth]{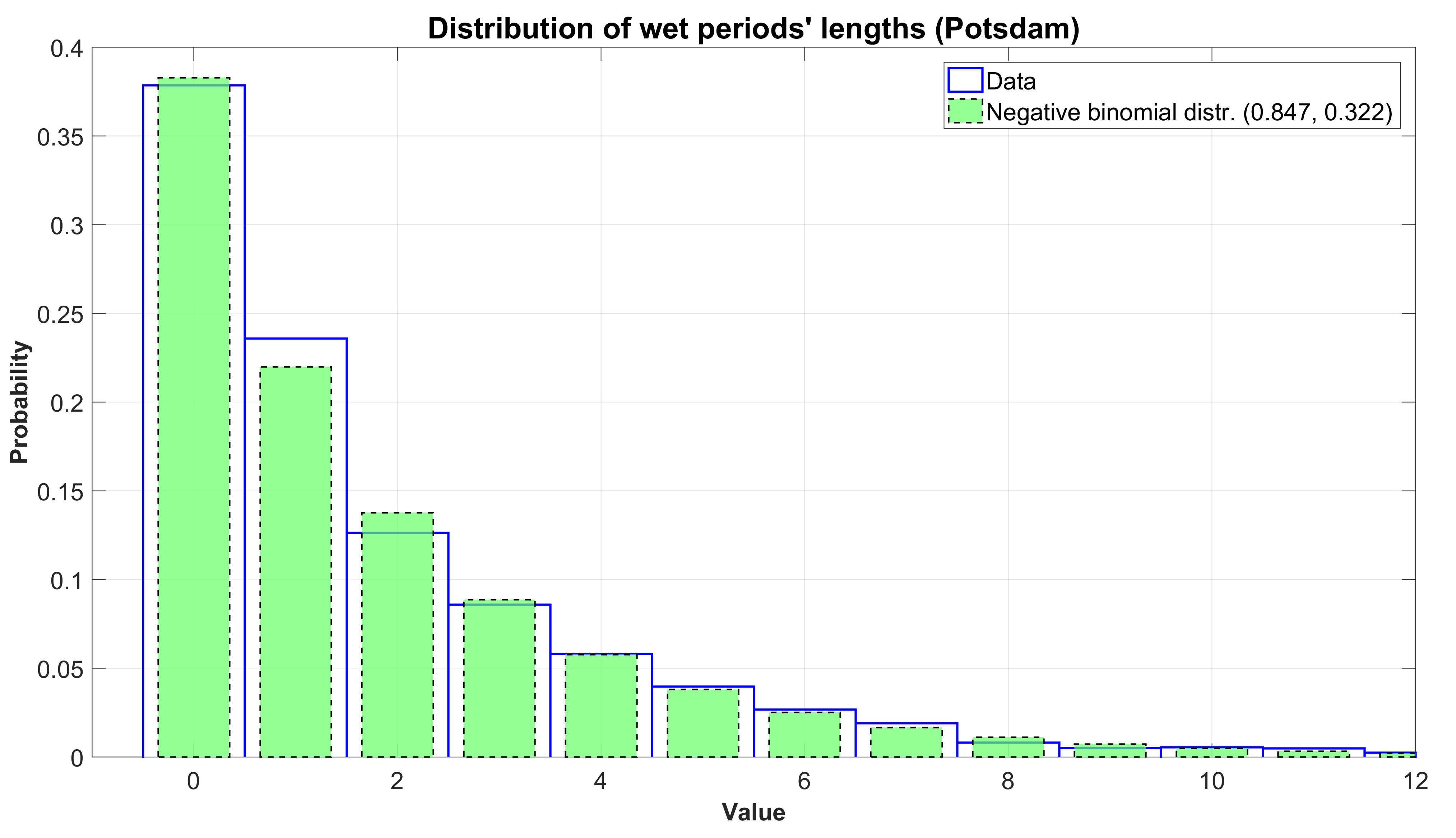} \\ b)}
\end{minipage}
\label{WetHist} \caption{The histogram of durations of wet periods
in Elista (a) and Potsdam (b) and the fitted negative binomial
distribution.}
\end{figure}

It is worth noting that at the same time the statistical analysis
convincingly suggests the Pareto-type model for the distribution of
daily precipitation volumes, see Figures~4a and 4b. For comparison,
on these figures there are also presented the graphs of the best
gamma-densities which, nevertheless, fit the histograms in a
noticeably worse way than the Pareto distributions.

\begin{figure}[h]
\begin{minipage}[h]{0.49\textwidth}
\center{\includegraphics[width=\textwidth,
height=0.6\textwidth]{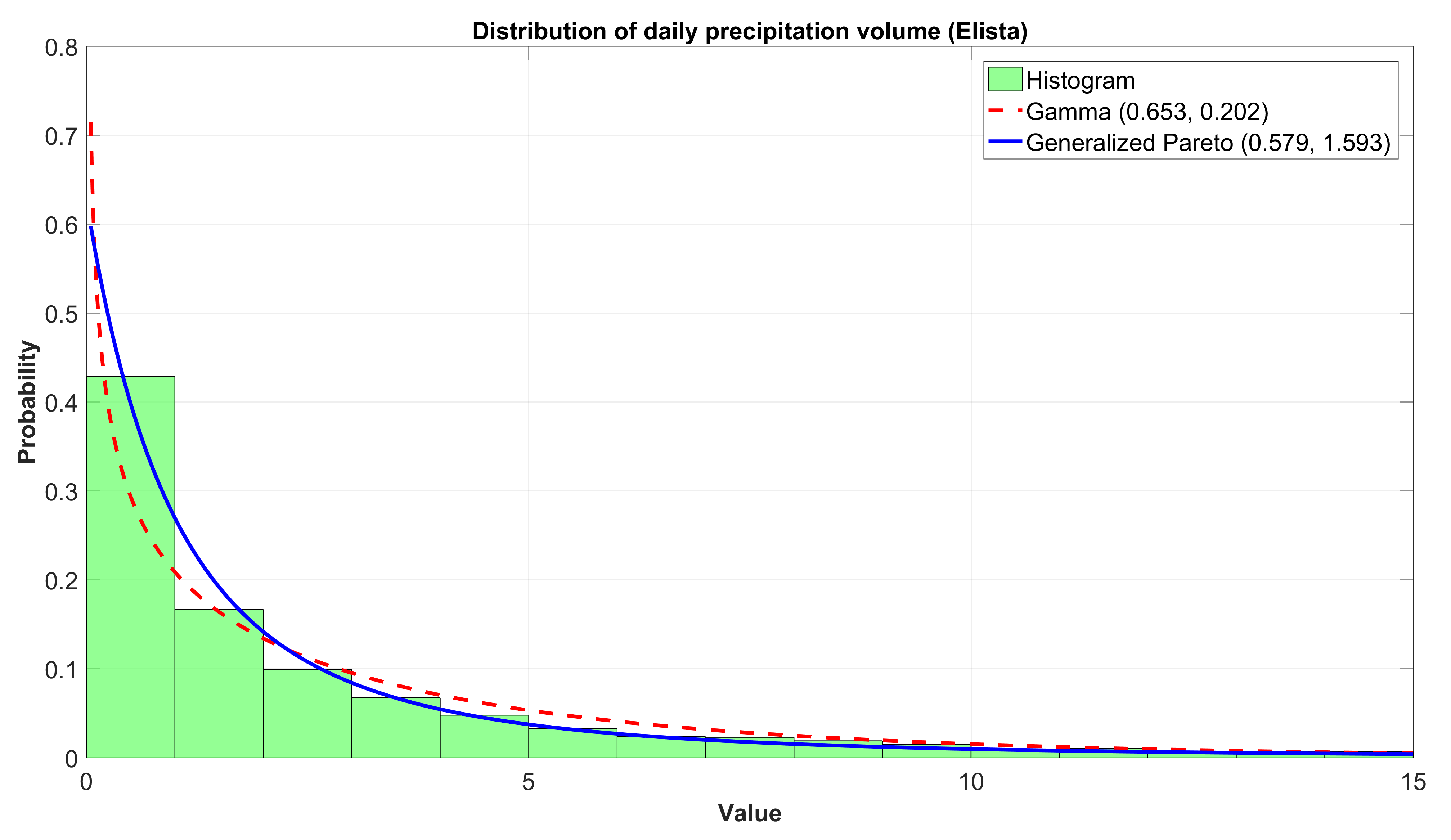}
\\a)}
\end{minipage}
\hfill
\begin{minipage}[h]{0.49\textwidth}
\center{
\includegraphics[width=\textwidth,
height=0.6\textwidth]{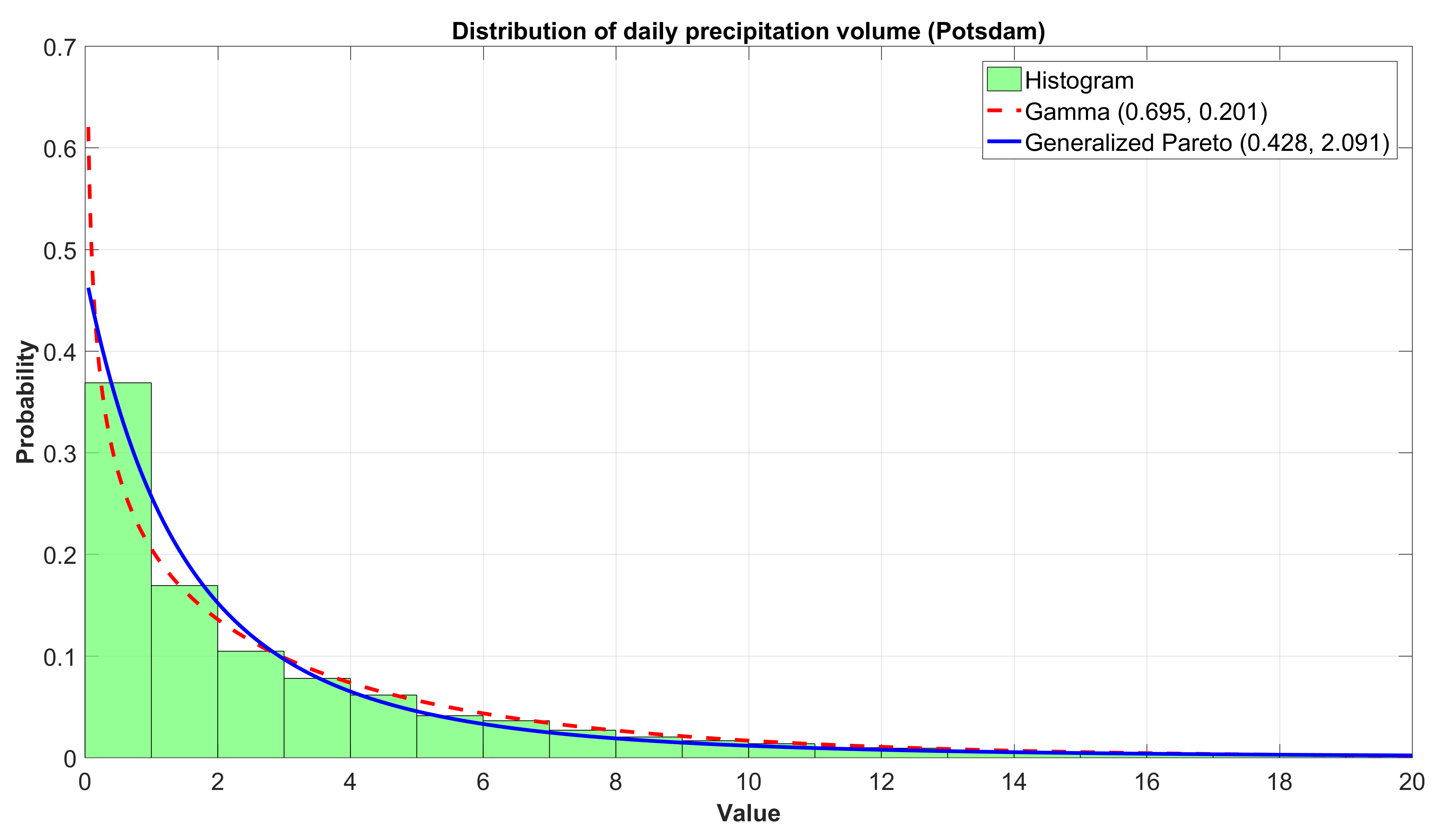} \\ b)}
\end{minipage}
\label{WetHist} \caption{The histogram of daily precipitation
volumes in Elista (a) and Potsdam (b) and the fitted Pareto and
gamma distributions.}
\end{figure}

In the same paper a schematic attempt was undertaken to explain this
phenomenon by the fact that negative binomial distributions can be
represented as mixed Poisson laws with mixing gamma-distributions.
As is known, the Poisson distribution is the best model for the
discrete stochastic chaos~\cite{Kingman1993} by virtue of the
universal principle of non-decrease of entropy in closed systems
(see, e. g., \cite{GnedenkoKorolev1996, KorolevBeningShorgin2011})
and the mixing distribution accumulates the statistical regularities
in the influence of stochastic factors that can be assumed exogenous
with respect to the local system under consideration.

In the paper \cite{Korolev2017} this explanation of the adequacy of
the negative binomial model was concretized. For this purpose, the
concept of a mixed geometric distribution introduced
in~\cite{Korolev2016TVP} (also see~\cite{KorolevPoisson,
Korolev2016}) was used. In~\cite{Korolev2017} it was demonstrated
that any negative binomial distribution with shape parameter no
greater than one is a mixed geometric distribution (this result is
reproduced below as Theorem 1). Thereby, a ``discrete'' analog of a
theorem due to L.~Gleser~\cite{Gleser1989} was proved. Gleser's
theorem establishes that a gamma distribution with shape parameter
no greater than one can be represented as a mixed exponential
distribution.

The representation of a negative binomial distribution as a mixed
geometric law can be interpreted in terms of the Bernoulli trials as
follows. First, as a result of some ``preliminary'' experiment the
value of some random variables (r.v:s) taking values in $[0,1]$ is
determined which is then used as the probability of success in the
sequence of Bernoulli trials in which the original ``unconditional''
r.v. with the negative binomial distribution is nothing else than
the ``conditionally'' geometrically distributed r.v. having the
sense of the number of trials up to the first failure. This makes it
possible to assume that the sequence of wet/dry days is not
independent, but is conditionally independent and the random
probability of success is determined by some outer stochastic
factors. As such, we can consider the seasonality or the type of the
cause of a rainy period.

The negative binomial model for the distribution of the duration of
wet periods makes it possible to obtain asymptotic approximations
for important characteristics of precipitation such as the
distribution of the total precipitation volume per wet period and
the distribution of the maximum daily precipitation volume within a
wet period. The first of these approximations was proposed
in~\cite{Korolev2017}, where an analog of the law of large numbers
for negative binomial random sums was presented stating that the
limit distribution for these sums is the gamma distribution.

The construction of the second approximation is the target of the
present paper.

The paper is organized as follows. Definitions and notation are
introduced in Section~2 which also contains some preliminary results
providing some theoretical grounds for the negative binomial model
of the probability distribution of the duration of wet periods. Main
results are presented and proved in Section 3 where the asymptotic
approximation is proposed for the distribution of the maximum daily
precipitation volume within a wet period. Some analytic properties
of the obtained limit distribution are described. In particular, it
is demonstrated that under certain conditions the limit distribution
is mixed exponential and hence, is infinitely divisible. It is shown
that under the same conditions the limit distribution can be
represented as a scale mixture of stable or Weibull or Pareto or
folded normal laws. The corresponding product representations for
the limit random variable can be used for its computer simulation.
Several methods for the statistical estimation of the parameters of
this distribution are proposed in Section 4. Section 5 contains the
results of fitting the distribution proposed in Section 3 to real
data by the methods described in Section 4.

\section{Preliminaries}

Although the main objects of our interest are the probability
distributions, for convenience and brevity in what follows we will
expound our results in terms of r.v:s with the corresponding
distributions assuming that all the r.v:s under consideration are
defined on one and the same probability space
$(\Omega,\,\mathfrak{F},\,{\sf P})$.

In the paper, conventional notation is used. The symbols $\eqd$ and
$\Longrightarrow$ denote the coincidence of distributions and
convergence in distribution, respectively. The integer and
fractional parts of a number $z$ will be respectively denoted $[z]$
and $\{z\}$.

A r.v. having the gamma distribution with shape parameter $r>0$ and
scale parameter $\lambda>0$ will be denoted $G_{r,\lambda}$,
$$
{\sf P}(G_{r,\lambda}<x)=\int_{0}^{x}g(z;r,\lambda)dz,\ \
\text{with}\ \
g(x;r,\lambda)=\frac{\lambda^r}{\Gamma(r)}x^{r-1}e^{-\lambda x},\
x\ge0,
$$
where $\Gamma(r)$ is Euler's gamma-function,
$\Gamma(r)=\int_{0}^{\infty}x^{r-1}e^{-x}dx$, $r>0$.

In these notation, obviously, $G_{1,1}$ is a r.v. with the standard
exponential distribution: ${\sf P}(G_{1,1}<x)=\big[1-e^{-x}\big]{\bf
1}(x\ge0)$ (here and in what follows ${\bf 1}(A)$ is the indicator
function of a set $A$).

The gamma distribution is a particular representative of the class
of {\it generalized gamma distributions} (GG-distributions), which
were first described in \cite{Stacy1962} as a special family of
lifetime distributions containing both gamma distributions and
Weibull distributions. A GG-distribution is the absolutely
continuous distribution defined by the density
$$
g^*(x;r,\gamma,\lambda)=\frac{|\gamma|\lambda^r}{\Gamma(r)}x^{\gamma
r-1}e^{-\lambda x^{\gamma}},\ \ \ \ x\ge0,
$$
with $\gamma\in\mathbb{R}$, $\lambda>0$, $r>0$.

The properties of GG-distributions are described in \cite{Stacy1962,
KorolevZaks2013}. A r.v. with the density $g^*(x;r,\gamma,\lambda)$
will be denoted $G^*_{r,\gamma,\lambda}$. It can be easily made sure
that
\begin{equation}\label{GG}
G^*_{r,\gamma,\lambda}\eqd G_{r,\lambda}^{1/\gamma}.
\end{equation}
For a r.v. with the Weibull distribution, a particular case of
GG-distributions corresponding to the density $g^*(x;1,\gamma,1)$
and the distribution function (d.f.)
$\big[1-e^{-x^{\gamma}}\big]{\bf 1}(x\ge0)$, we will use a special
notation $W_{\gamma}$. Thus, $G_{1,1}\eqd W_1$. It is easy to see
that
\begin{equation}\label{Weibull}
W_1^{1/\gamma}\eqd W_{\gamma}.
\end{equation}
A r.v. with the standard normal d.f. $\Phi(x)$ will be denoted $X$,
$$
{\sf
P}(X<x)=\Phi(x)=\frac{1}{\sqrt{2\pi}}\int_{-\infty}^{x}e^{-z^2/2}dz,\
\ \ \ x\in\mathbb{R}.
$$
The distribution of the r.v. $|X|$, i. e. ${\sf
P}(|X|<x)=2\Phi(x)-1$, $x\ge0$, is called half-normal or folded
normal.

The d.f. and the density of a strictly stable distribution with the
characteristic exponent $\alpha$ and shape parameter $\theta$
defined by the characteristic function (ch.f.)
$$
\mathfrak{f}_{\alpha,\theta}(t)=\exp\big\{-|t|^{\alpha}\exp\{-{\textstyle\frac12}i\pi\theta\alpha\mathrm{sign}t\}\big\},\
\ \ \ t\in\r,
$$
wheter $0<\alpha\le2$, $|\theta|\le\min\{1,\frac{2}{\alpha}-1\}$,
will be respectively denoted $F_{\alpha,\theta}(x)$ and
$f_{\alpha,\theta}(x)$ (see, e. g., \cite{Zolotarev1983}). A r.v.
with the d.f. $F_{\alpha,\theta}(x)$ will be denoted
$S_{\alpha,\theta}$.

To symmetric strictly stable distributions there correspond the
value $\theta=0$ and the ch.f.
$\mathfrak{f}_{\alpha,0}(t)=e^{-|t|^{\alpha}}$, $t\in\r$.
To one-sided strictly stable distributions concentrated on the
nonnegative halfline there correspond the values $\theta=1$ and
$0<\alpha\le1$. The pairs $\alpha=1$, $\theta=\pm1$ correspond to
the distributions degenerate in $\pm1$, respectively. All the rest
strictly stable distributions are absolutely continuous. Stable
densities cannot be explicitly represented via elementary functions
with four exceptions: the normal distribution ($\alpha=2$,
$\theta=0$), the Cauchy distribution ($\alpha=1$, $\theta=0$), the
L{\'e}vy distribution ($\alpha=\frac12$, $\theta=1$) and the
distribution symmetric to the L{\'e}vy law ($\alpha=\frac12$,
$\theta=-1$).

It is well known that if $0<\alpha<2$, then ${\sf
E}|S_{\alpha,\theta}|^{\beta}<\infty$ for any $\beta\in(0,\alpha)$,
but the moments of the r.v. $S_{\alpha,\theta}$ of orders
$\beta\ge\alpha$ do not exist (see, e. g., \cite{Zolotarev1983}).
Despite the absence of explicit expressions for the densities of
stable distributions in terms of elementary functions, it can be
shown \cite{KorolevWeibull2016} that
\begin{equation}
{\sf
E}S_{\alpha,1}^{\beta}=\frac{\Gamma\big(1-\frac{\beta}{\alpha}\big)}{\Gamma(1-\beta)}
\label{stabmom}
\end{equation}
for $0<\beta<\alpha\le 1$.

In \cite{KotzOstrovskii1996, KorolevZeifman2016a,
KorolevZeifman2016b} it was proved that if $\alpha\in(0,1)$ and the
i.i.d. r.v:s $S_{\alpha,1}$ and $S'_{\alpha,1}$ have the same
strictly stable distribution, then the density $v_{\alpha}(x)$ of
the r.v. $R_{\alpha}=S_{\alpha,1}/S'_{\alpha,1}$ has the form
\begin{equation}
v_{\alpha}(x)=\frac{\sin(\pi\alpha)x^{\alpha-1}}{\pi[1+x^{2\alpha}+2x^{\alpha}\cos(\pi\alpha)]},\
\ \ x>0.
\label{Rdensity}
\end{equation}

A r.v. $N_{r,p}$ is said to have the {\it negative binomial
distribution} with parameters $r>0$ (``shape'') and $p\in(0,1)$
(``success probability''), if
$$
{\sf P}(N_{r,p}=k)=\frac{\Gamma(r+k)}{k!\Gamma(r)}\cdot p^r(1-p)^k,\
\ \ \ k=0,1,2,...
$$

A particular case of the negative binomial distribution
corresponding to the value $r=1$ is the {\it geometric
distribution}. Let $p\in(0,1)$ and let $N_{1,p}$ be the r.v. having
the geometric distribution with parameter $p\,$:
$$
{\sf P}(N_{1,p}=k)=p(1-p)^{k},\ \ \ \ k=0,1,2,...
$$
This means that for any $m\in\mathbb{N}$
$$
{\sf P}(N_{1,p}\ge
m)=\sum\nolimits_{k=m}^{\infty}p(1-p)^{k}=(1-p)^{m}.
$$

Let $Y$ be a r.v. taking values in the interval $(0,1)$. Moreover,
let for all $p\in(0,1)$ the r.v. $Y$ and the geometrically
distributed r.v. $N_{1,p}$ be independent. Let $V=N_{1,Y}$, that is,
$V(\omega)=N_{1,Y(\omega)}(\omega)$ for any $\omega\in\Omega$. The
distribution
$$
{\sf P}(V\ge m)=\int_{0}^{1}(1-y)^{m}d{\sf P}(Y<y), \ \ \
m\in\mathbb{N},
$$
of the r.v. $V$ will be called {\it $Y$-mixed geometric}
\cite{Korolev2016TVP}.

It is well known that the negative binomial distribution is a mixed
Poisson distribution with the gamma mixing distribution
\cite{GreenwoodYule1920} (also see \cite{KorolevBeningShorgin2011}):
for any $r>0$, $p\in(0,1)$ and $k\in\{0\}\bigcup\mathbb{N}$ we have
\begin{equation}
\frac{\Gamma(r+k)}{k!\Gamma(r)}\cdot
p^r(1-p)^k=\frac{1}{k!}\int_{0}^{\infty}e^{-z}z^kg(z;r,\mu)dz,\label{NBMixt}
\end{equation}
where $\mu=p/(1-p)$.

Based on representation \eqref{NBMixt}, in \cite{Korolev2017} it was
proved that any negative binomial distribution with the shape
parameter no greater than one is a mixed geometric distribution.
Namely, the following statement was proved that gives an analytic
explanation of the validity of the negative binomial model for the
duration of wet periods measured in days (see the Introduction).

\smallskip

{\sc Theorem 1} \cite{Korolev2017}. {\it The negative binomial
distribution with parameters $r\in(0,1)$ and $p\in(0,1)$ is a mixed
geometric distribution$:$ for any $k\in\{0\}\bigcup\mathbb{N}$}
$$
\frac{\Gamma(r+k)}{k!\Gamma(r)}\cdot
p^r(1-p)^k=\int_{\mu}^{\infty}\Big(\frac{z}{z+1}\Big)\Big(1-\frac{z}{z+1}\Big)^kp(z;r,\mu)dz=\int_{p}^{1}y(1-y)^kh(y;r,p)dy,
$$
{\it where $\mu=p/(1-p)$ and the probability densities $p(z;r,\mu)$
and $h(y;r,p)$ have the forms
$$
p(z;r,\mu)=\frac{\mu^r}{\Gamma(1-r)\Gamma(r)}\cdot\frac{\mathbf{1}(z\ge\mu)}{(z-\mu)^rz},
$$
$$
h(y;r,p)=\frac{p^r}{\Gamma(1-r)\Gamma(r)}\cdot\frac{(1-y)^{r-1}\mathbf{1}(p<y<1)}{y(y-p)^r}.
$$

Furthermore, if $G_{r,\,1}$ and $G_{1-r,\,1}$ are independent
gamma-distributed r.v:s, $\mu>0$, $p\in(0,1)$, then the density
$p(z;r,\mu)$ corresponds to the r.v.
\begin{equation}
Z_{r,\mu}=\frac{\mu(G_{r,\,1}+G_{1-r,\,1})}{G_{r,\,1}}
\label{Zdef}
\end{equation}
and the density $h(y;r,p)$ corresponds to the r.v.
$$
Y_{r,p}=\frac{p(G_{r,\,1}+G_{1-r,\,1})}{G_{r,\,1}+pG_{1-r,\,1}}.
$$
}

\smallskip

Let $P(t)$, $t\ge0$, be the standard Poisson process (homogeneous
Poisson process with unit intensity). Then
distribution~\eqref{NBMixt} corresponds to the r.v.
$N_{r,p}=P(G_{r,p/(1-p)})$, where the r.v. $G_{r,p/(1-p)}$ is
independent of the process $P(t)$.

\section{The probability distribution of extremal precipitation}

In this section we will deduce the probability distribution of
extremal daily precipitation within a wet period.

Let $r>0$, $\lambda>0$, $q\in(0,1)$, $n\in\mathbb{N}$,
$p_n=\min\{q,\,\lambda/n\}$. It is easy to make sure that
\begin{equation}
n^{-1}G_{r,p_n/(1-p_n)}\Longrightarrow G_{r,\lambda}\label{2}
\end{equation}
as $n\to\infty$.

\smallskip

{\sc Lemma 1.} {\it Let $\Lambda_1,\Lambda_2,\ldots$ be a sequence
of positive r.v$:$s such that for any $n\in\mathbb{N}$ the r.v.
$\Lambda_n$ is independent of the Poisson process $P(t)$, $t\ge0$.
The convergence
$$
n^{-1}P(\Lambda_n)\Longrightarrow \Lambda
$$
as $n\to\infty$ to some nonnegative r.v. $\Lambda$ takes place if
and only if
\begin{equation}
n^{-1}\Lambda_n\Longrightarrow \Lambda \label{3}
\end{equation}
as $n\to\infty$.}

\smallskip

{\sc Proof}. This statement is a particular case of Lemma 2
in~\cite{Korolev1998} (also see Theorem 7.9.1 in
\cite{KorolevBeningShorgin2011}).

\smallskip

Consider a sequence of independent identically distributed (i.i.d.)
r.v:s $X_1,X_2,\ldots$. Let $N_1,N_2,\ldots$ be a sequence of
natural-valued r.v:s such that for each $n\in\mathbb{N}$ the r.v.
$N_n$ is independent of the sequence $X_1,X_2,\ldots$. Denote
$M_n=\max\{X_1,\ldots,X_{N_n}\}$.

Let $F(x)$ be a d.f., $a\in\mathbb{R}$. Denote
$\mathrm{rext}(F)=\sup\{x:\,F(x)<1\}$, $F^{-1}(a)=\inf\{x:\,F(x)\ge
a\}$.

\smallskip

{\sc Lemma 2.} {\it Let $\Lambda_1,\Lambda_2,\ldots$ be a sequence
of positive r.v$:$s such that for each $n\in\mathbb{N}$ the r.v.
$\Lambda_n$ is independent of the Poisson process $P(t)$, $t\ge0$.
Let $N_n=P(\Lambda_n)$. Assume that there exists a nonnegative r.v.
$\Lambda$ such that convergence~{\rm \eqref{3}} takes place. Let
$X_1,X_2,\ldots$ be i.i.d. r.v$:$s with a common d.f. $F(x)$. Assume
also that $\mathrm{rext}(F)=\infty$ and there exists a number
$\gamma>0$ such that for each $x>0$
\begin{equation}
\lim_{y\to\infty}\frac{1-F(xy)}{1-F(y)}=x^{-\gamma}.\label{4}
\end{equation}
Then}
$$
\lim_{n\to\infty}\sup_{x\ge 0}\bigg|{\sf
P}\bigg(\frac{M_n}{F^{-1}(1-\frac{1}{n})}<x\bigg)-\int_{0}^{\infty}e^{-zx^{-\gamma}}d{\sf
P}(\Lambda<z)\bigg|=0.
$$

\smallskip

{\sc Proof.} This statement is a particular case of Theorem 3.1 in
\cite{KorolevSokolov2008}.

\smallskip

{\sc Theorem 2.} {\it Let $n\in\mathbb{N}$, $\lambda>0$, $q\in(0,1)$
and let $N_{r,p_n}$ be a r.v. with the negative binomial
distribution with parameters $r>0$ and $p_n=\min\{q,\lambda/n\}$.
Let $X_1,X_2,\ldots$ be i.i.d. r.v$:$s with a common d.f. $F(x)$.
Assume that $\mathrm{rext}(F)=\infty$ and there exists a number
$\gamma>0$ such that relation~{\rm \eqref{4}} holds for any $x>0$.
Then
$$
\lim_{n\to\infty}\sup_{x\ge 0}\bigg|{\sf
P}\bigg(\frac{\max\{X_1,\ldots,X_{N_{r,p_n}}\}}{F^{-1}(1-\frac{1}{n})}<x\bigg)-F(x;
r,\lambda,\gamma)\bigg|=0,
$$
where}
$$
F(x; r,\lambda,\gamma)=\bigg(\frac{\lambda x^{\gamma}}{1+\lambda
x^{\gamma}}\bigg)^r,\ \ \ x\ge0.
$$

\smallskip

{\sc Proof.} From~\eqref{NBMixt} it follows that $N_{r,p_n}\eqd
P(G_{r,p_n/(1-p_n)}$. Therefore, from~\eqref{2}, Lemma 1 with
$\Lambda_n=G_{r,p_n/(1-p_n)}$ and Lemma 2 with the account of the
absolute continuity of the limit distribution it immediately follows
that
$$
\lim_{n\to\infty}\sup_{x\ge 0}\bigg|{\sf
P}\bigg(\frac{\max\{X_1,\ldots,X_{N_{r,p_n}}\}}{F^{-1}(1-\frac{1}{n})}<x\bigg)-
\frac{\lambda^r}{\Gamma(r)}\int_{0}^{\infty}e^{-z(\lambda+x^{-\gamma})}z^{r-1}dz\bigg|=0.
$$
By elementary calculation it can be made sure that
$$
\frac{\lambda^r}{\Gamma(r)}\int_{0}^{\infty}e^{-z(\lambda+x^{-\gamma})}z^{r-1}dz=
\frac{\lambda^r}{\Gamma(r)(\lambda+x^{-\gamma})^r}\int_{0}^{\infty}e^{-z}z^{r-1}dz=\bigg(\frac{\lambda
x^{\gamma}}{1+\lambda x^{\gamma}}\bigg)^r.
$$
The theorem is proved.

\smallskip

It is worth noting that the limit distribution with the power-type
decrease of the tail was obtained in Theorem~2 as a scale mixture of
the Fr{\'e}chet distribution (the type II extreme value
distribution) in which the mixing law is the gamma distribution.

Since the Fr{\'e}chet d.f. $e^{-x^{-\gamma}}$ with $\gamma>0$
corresponds to the r.v. $W_{\gamma}^{-1}$, it is easy to make sure
that the d.f. $F(x; r,\lambda,\gamma)$ corresponds to the r.v.
$M_{r,\gamma,\lambda}\equiv
G_{r,\lambda}^{1/\gamma}W_{\gamma}^{-1}$, where the multipliers on
the right-hand side are independent. From~\eqref{GG}
and~\eqref{Weibull} it follows that
\begin{equation}\label{M}
M_{r,\gamma,\lambda}\eqd\Big(\frac{G_{r,\lambda}}{W_1}\Big)^{1/\gamma}
\eqd\frac{G^*_{r,\gamma,\lambda}}{W_{\gamma}}
\end{equation}
where in each term the multipliers are independent. Consider the
r.v. $G_{r,\lambda}/W_1$ in \eqref{M} in more detail. We have
$$
\frac{G_{r,\lambda}}{W_1}\eqd\frac{G_{r,\lambda}}{G_{1,1}}\eqd\frac{G_{r,1}}{\lambda
G_{1,1}}\eqd\frac{Q_{r,1}}{\lambda r},
$$
where $Q_{r,1}$ is the r.v. having the Snedecor--Fisher distribution
with parameters $r,\,1$ (`degrees of freedom') defined by the
Lebesgue density
$$
f_{r,1}(x)=\frac{r^{r+1}x^{r-1}}{(1+rx)^{r+1}},\ \ \ x\ge0,
$$
(see, e. g., \cite{Bolshev}, Section 27).

So,
\begin{equation}\label{MQ}
M_{r,\gamma,\lambda}\eqd\Big(\frac{Q_{r,1}}{\lambda
r}\Big)^{1/\gamma},
\end{equation}
and the statement of theorem 2 can be re-formulated as
\begin{equation}
\label{Mdef}
\frac{\max\{X_1,\ldots,X_{N_{r,p_n}}\}}{F^{-1}(1-\frac{1}{n})}\Longrightarrow
M_{r,\gamma,\lambda}\equiv
\frac{G_{r,\lambda}^{1/\gamma}}{W_{\gamma}}\eqd
\Big(\frac{Q_{r,1}}{\lambda r}\Big)^{1/\gamma}\ \ \ \ (n\to\infty).
\end{equation}

The density of the limit distribution $F(x;r,\gamma,\lambda)$ of the
extreme daily precipitation within a wet period has the form
\begin{equation}
p(x;r,\gamma,\lambda)=\frac{r\gamma\lambda^rx^{\gamma
r-1}}{(1+\lambda x^{\gamma})^{r+1}}=\frac{\gamma
r\lambda^r}{x^{1+\gamma}(\lambda+x^{-\gamma})^{r+1}},\ \ \
x>0.\label{ExtrPDF}
\end{equation}

It is easy to see that $p(x;r,\gamma,\lambda)=O(x^{-1-\gamma})$ as
$x\to\infty$. Therefore ${\sf
E}M_{r,\gamma,\lambda}^{\delta}<\infty$ only if $\delta<\gamma$.
Moreover, from~\eqref{Mdef} it is possible to deduce explicit
expressions for the moments of the r.v. $M_{r,\gamma,\lambda}$.

\smallskip

{\sc Theorem 3.} {\it Let $0<\delta<\gamma<\infty$. Then}
$$
{\sf
E}M_{r,\gamma,\lambda}^{\delta}=
\frac{\Gamma\big(r+\frac{\delta}{\gamma}\big)\Gamma\big(1-\frac{\delta}{\gamma}\big)}{\lambda^{\delta/\gamma}\Gamma(r)}.
$$

\smallskip

{\sc Proof}. From \eqref{Mdef} it follows that
\begin{equation}
\label{Mmoments} {\sf E}M_{r,\gamma,\lambda}^{\delta}={\sf
E}G_{r,\lambda}^{\delta/\gamma}\cdot{\sf E}W_1^{-\delta/\gamma}.
\end{equation}
It is easy to verify that
\begin{equation}
\label{moments} {\sf
E}G_{r,\lambda}^{\delta/\gamma}=\frac{\Gamma\big(r+\frac{\delta}{\gamma}\big)}{\lambda^{\delta/\gamma}\Gamma(r)},\
\ \ {\sf
E}W_1^{-\delta/\gamma}=\Gamma\big(1-{\textstyle\frac{\delta}{\gamma}}\big).
\end{equation}
Hence follows the desired result.

\smallskip

To analyze the properties of the limit distribution in theorem 2
more thoroughly we will require some additional auxiliary results.

\smallskip

{\sc Lemma 3} \cite{KorolevWeibull2016}. {\it Let $\gamma\in(0,1]$.
Then
$$
W_{\gamma}\eqd \frac{W_1}{S_{\gamma,1}}
$$
with the r.v:s on the right-hand side being independent.}

\smallskip

{\sc Lemma 4} \cite{Korolev2017}. {\it Let $r\in(0,1]$,
$\gamma\in(0,1]$, $\lambda>0$. Then
$$
G_{r,\lambda}^{1/\gamma}\eqd
G^*_{r,\gamma,\lambda}\eqd\frac{W_{\gamma}}{Z_{r,\lambda}^{1/\gamma}}\eqd
\frac{W_1}{S_{\gamma,1}Z_{r,\lambda}^{1/\gamma}},
$$
where the r.v. $Z_{r,\lambda}$ was defined in \eqref{Zdef} and all
the involved r.v$:$s are independent.}

\smallskip

{\sc Theorem 4}. {\it Let $r\in(0,1]$, $\gamma\in(0,1]$,
$\lambda>0$. Then the following product representations are valid$:$
\begin{equation}\label{T3_1}
M_{r,\gamma,\lambda}\eqd
\frac{G_{r,\lambda}^{1/\gamma}S_{\gamma,1}}{W_1},
\end{equation}
\begin{equation}\label{T3_2}
M_{r,\gamma,\lambda}\eqd
\frac{W_{\gamma}}{W'_{\gamma}}\cdot\frac{1}{Z_{r,\lambda}^{1/\gamma}}\eqd
W_1\cdot\frac{R_{\gamma}}{W'_1Z_{r,\lambda}^{1/\gamma}}\eqd
\frac{\Pi R_{\gamma}}{Z_{r,\lambda}^{1/\gamma}}\eqd
\frac{|X|\sqrt{2W_1}R_{\gamma}}{W'_1Z_{r,\lambda}^{1/\gamma}},
\end{equation}
where $W_{\gamma}\eqd W'_{\gamma}$, $W_1\eqd W'_1$, the r.v.
$R_{\gamma}$ has the density {\rm\eqref{Rdensity}}, the r.v. $\Pi$
has the Pareto distribution$:$ ${\sf P}(\Pi>x)=(x+1)^{-1}$, $x\ge0$,
and in each term the involved r.v$:$s are independent.}

\smallskip

{\sc Proof}. Relation \eqref{T3_1} follows from \eqref{Mdef} and
Lemma 3, relation \eqref{T3_2} follows from \eqref{Mdef} and Lemma 4
with the account of the representation $W_1\eqd |X|\sqrt{2W_1}$, the
proof of which can be found in, say, \cite{KorolevWeibull2016}.

\smallskip

With the account of the relation $R_{\gamma}\eqd R_{\gamma}^{-1}$,
from~\eqref{T3_2} we obtain the following statement.

\smallskip

{\sc Corollary 1.} {\it Let $r\in(0,1]$, $\gamma\in(0,1]$,
$\lambda>0$. Then the d.f. $F(x;r,\gamma,\lambda)$ is mixed
exponential$:$
$$
1-F(x;r,\gamma,\lambda)=\int_{0}^{\infty}e^{-ux}dA(u),\ \ \ x\ge0,
$$
where
$$
A(u)={\sf P}\big(W_1R_{\gamma}Z_{r,\lambda}^{1/\gamma}<u\big),\ \ \
u\ge0,
$$
and all the involved r.v$:$s are independent}.

\smallskip

{\sc Remark 1.} It is worth noting that the mixing distribution in
Corollary 1 is mixed exponential itself.

\smallskip

{\sc Corollary 2.} {\it Let $r\in(0,1]$, $\gamma\in(0,1]$,
$\lambda>0$. Then the d.f. $F(x;r,\gamma,\lambda)$ is infinitely
divisible.}

\smallskip

{\sc Proof.} This statement immediately follows from Corollary 1 and
the result of Goldie \cite{Goldie1967} stating that the product of
two independent non-negative random variables is infinitely
divisible, if one of the two is exponentially distributed.

\smallskip

Theorem 3 states that the limit distribution in Theorem 2 can be
represented as a scale mixture of exponential or stable or Weibull
or Pareto or folded normal laws. The corresponding product
representations for the r.v. $M_{r,\gamma,\lambda}$ can be used for
its computer simulation.

In practice, the asymptotic approximation $F(x; r,\lambda,\gamma)$
for the distribution of the extreme daily precipitation within a wet
period proposed by Theorem~2 is adequate, if the ``success
probability'' is small enough, that is, if on the average the wet
periods are long enough.

\section{Estimation of the parameters $r$, $\lambda$ and $\gamma$}

From~\eqref{ExtrPDF} it can be seen that the realization of the
maximum likelihood method for the estimation of the parameters $r$,
$\lambda$ and $\gamma$ inevitably assumes the necessity of numerical
solution of a system of transcendental equations by iterative
procedures without any guarantee that the resulting maximum is
global. The closeness of the initial approximation to the true
maximum likelihood point in the three-dimensional parameter set
might give a hope that the terminal extreme point found by the
numerical algorithm is global.

For rough estimation of the parameters, the following considerably
simpler method can be used. The resulting rough estimates can be
used as a starting point for the `full' maximum likelihood algorithm
mentioned above in order to ensure the closeness of the initial
approximation to the true solution. The rough method is based on
that the quantiles of the d.f. $F(x; r,\lambda,\gamma)$ can be
written out explicitly. Namely, the quantile
$x(\epsilon;r,\lambda,\gamma)$ of the d.f. $F(x; r,\lambda,\gamma)$
of order $\epsilon\in(0,1)$, that is, the solution of the equation
$F(x; r,\lambda,\gamma)=\epsilon$ with respect to $x$, obviously has
the form
$$
x(\epsilon;r,\lambda,\gamma)=\bigg(\frac{\epsilon^{1/r}}{\lambda-\lambda\epsilon^{1/r}}\bigg)^{1/\gamma}.
$$
Let at our disposal there be observations $\{X_{i,j}\}$,
$i=1,\ldots,m$, $j=1,\ldots,m_i$, where $i$ is the number of a wet
period (the number of a sequence of rainy days), $j$ is the number
of a day in the wet sequence, $m_i$ is the length of the $i$th wet
sequence (the number of rainy days in the $i$th wet period), $m$ is
the total number of wet sequences, $X_{i,j}$ is the precipitation
volume on the $j$th day of the $i$th wet sequence. Construct the
sample $X^*_1,\ldots,X^*_m$ as
\begin{equation}
X^*_k=\max\{X_{k,1},\ldots,X_{k,m_k}\},\ \ \ k=1,\ldots,m.\label{VarSample}
\end{equation}
Let $X^*_{(1)},\ldots,X^*_{(m)}$ be order statistics constructed
from the sample $X^*_1,\ldots,X^*_m$. Since we have three unknown
parameters $r$, $\lambda$ and $\gamma$, fix three numbers
$0<p_1<p_2<p_3<1$ and construct the system of equations
$$
{\displaystyle X^*_{([mp_k])}=\bigg(\frac{p_k^{1/r}}{\lambda-\lambda
p_k^{1/r}}\bigg)^{1/\gamma}},\ \ \ k=1,2,3
$$
(here the symbol $[a]$ denotes the integer part of a number $a$).

This system can be solved by standard techniques. For example,
first, the number $s\equiv\frac1r$ can be found numerically as the
solution of the equation
$$
Cs=\log\frac{1-p_3^s}{1-p_1^s}\log\frac{ X^*_{([mp_1])}}{
X^*_{([mp_2])}}-\log\frac{1-p_2^s}{1-p_1^s}\log\frac{
X^*_{([mp_1])}}{X^*_{([mp_3])}},
$$
where
$$
C=\log
\frac{X^*_{([mp_1])}}{X^*_{([mp_3])}}\log\frac{p_1}{p_2}-\log\frac{
X^*_{([mp_1])}}{X^*_{([mp_2])}}\log\frac{p_1}{p_3}.
$$
Having obtained the value of $s$, one can then find the values of
$\gamma$ and $\lambda$ explicitly:
\begin{equation}
\gamma=\frac{s(\log p_1-\log p_3)+\log(1-p_3^s)-\log(1-p_1^s)}{\log
X^*_{([mp_1])}-\log X^*_{([mp_3])}}, \label{gamma}
\end{equation}
\begin{equation}
\lambda=\frac{p_2^s}{(1-p_2^s)(X^*_{([mp_2])})^{\gamma}}.
\label{lambda}
\end{equation}
As $p_1$, $p_2$ and $p_3$ one may take, say, $p_k=\frac{k}{4}$,
$k=1,2,3$. Or it is possible to fix a $\tau\in(0,\frac14)$, set
$p_1=\tau$, $p_2=\frac12$, $p_3=1-\tau$ and choose a $\tau$ that
provides the best fit of the estimated model d.f. $F(x;
r,\lambda,\gamma)$ with the empirical d.f.

If the parameter $r$ is known (for example, it is estimated
beforehand while solving the problem of fitting the negative
binomial distribution to the empirical data on the lengths of wet
periods), then the parameters $\lambda$ and $\gamma$ can be
estimated directly by formulas~\eqref{gamma} and~\eqref{lambda}.

With known $r$, the more accurate estimates of the parameters
$\lambda$ and $\gamma$ can be also found by minimizing the
discrepancy between the empirical and model d.f:s by the least
squares techniques. Namely, from the Glivenko theorem it follows
that
\begin{equation}
\bigg(\frac{\lambda (X^*_{(i)})^{\gamma}}{1+\lambda
(X^*_{(i)})^{\gamma}}\bigg)^r\approx\frac{i}{m}.\label{7}
\end{equation}
As this is so, for every $i=1,\ldots,m-1$ the following implications
take place:
$$
\eqref{7}\Longleftrightarrow \Big\{\frac{\lambda
(X^*_{(i)})^{\gamma}}{1+\lambda
(X^*_{(i)})^{\gamma}}\approx\Big(\frac{i}{m}\Big)^{1/r}\Big\}\Longleftrightarrow\Big\{
\frac{1}{1+\lambda (X^*_{(i)})^{\gamma}}\approx
1-\Big(\frac{i}{m}\Big)^{1/r}\Big\}\Longleftrightarrow
$$
$$
\Longleftrightarrow\Big\{\lambda(X^*_{(i)})^{\gamma}\approx\frac{i^{1/r}}{m^{1/r}-i^{1/r}}\Big\}\Longleftrightarrow
\Big\{\log\lambda+\gamma\log
X^*_{(i)}\approx\log\frac{i^{1/r}}{m^{1/r}-i^{1/r}}\Big\}.
$$
Therefore the estimates $\widehat\lambda$ and $\widehat\gamma$ of
the parameters $\lambda$ and $\gamma$ can be found as the solution
of the least squares problem
$$
(\log\widehat\lambda,\,\widehat\gamma)=\mathrm{arg}\min_{\log\lambda,\gamma}\sum_{i=1}^{m-1}\big(\log\lambda+\gamma\log
X^*_{(i)}-c_i\big)^2,
$$
where
$$ c_i=\log\frac{i^{1/r}}{m^{1/r}-i^{1/r}}.
$$
This solution can be written out explicitly so that finally we have
\begin{equation}\label{gammaLS}
\widehat\gamma=\frac{(m-1)\sum_{i=1}^{m-1}c_i\log
X^*_{(i)}-\sum_{i=1}^{m-1}\log
X^*_{(i)}\sum_{i=1}^{m-1}c_i}{(m-1)\sum_{i=1}^{m-1}(\log
X^*_{(i)})^2-(\sum_{i=1}^{m-1}\log X^*_{(i)})^2},
\end{equation}
\begin{equation}\label{lambdaLS}
\widehat\lambda=\exp\Big\{\frac{1}{m-1}\Big(\sum\nolimits_{i=1}^{m-1}c_i-\widehat\gamma\sum\nolimits_{i=1}^{m-1}\log
X^*_{(i)}\Big)\Big\}.
\end{equation}

\section{The statistical analysis of real data}

In this section we present the results of statistical estimation of
the distribution of extremal daily precipitation within a wet period
by the methods described in the preceding section. As the data, we
use the observations of daily precipitation in Potsdam and Elista
from $1950$ to $2009$.

First of all, notice that the Pareto distribution of daily
precipitation volumes (see Figure 4) satisfies condition \eqref{4}.
Therefore, the asymptotic approximation provided by Theorem 2 can be
applied to the statistical analysis of the real data.

\begin{figure}[ptb]
\begin{minipage}[h]{0.44\linewidth}
\center{\includegraphics[width=1\linewidth, height=0.7\linewidth]{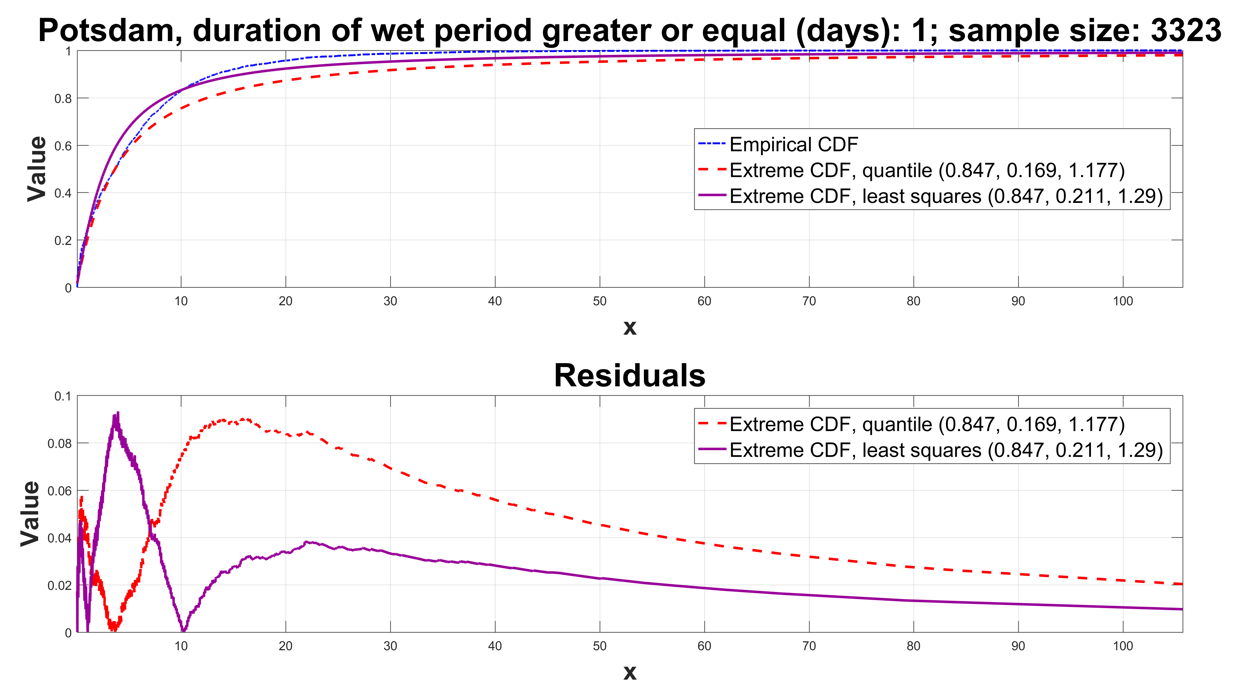}} a) \\
\end{minipage}
\hfill
\begin{minipage}[h]{0.44\linewidth}
\center{\includegraphics[width=1\linewidth,
height=0.7\linewidth]{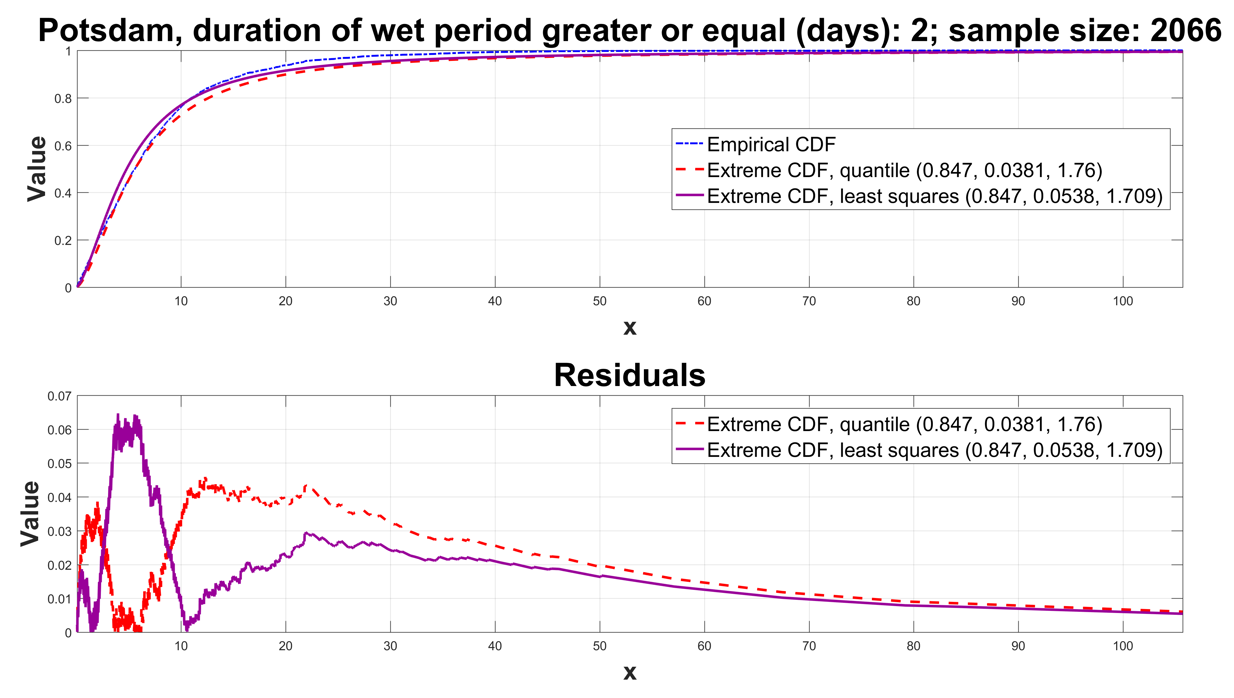}} b)\\
\end{minipage}
\vfill
\begin{minipage}[h]{0.44\linewidth}
\center{\includegraphics[width=1\linewidth, height=0.7\linewidth]{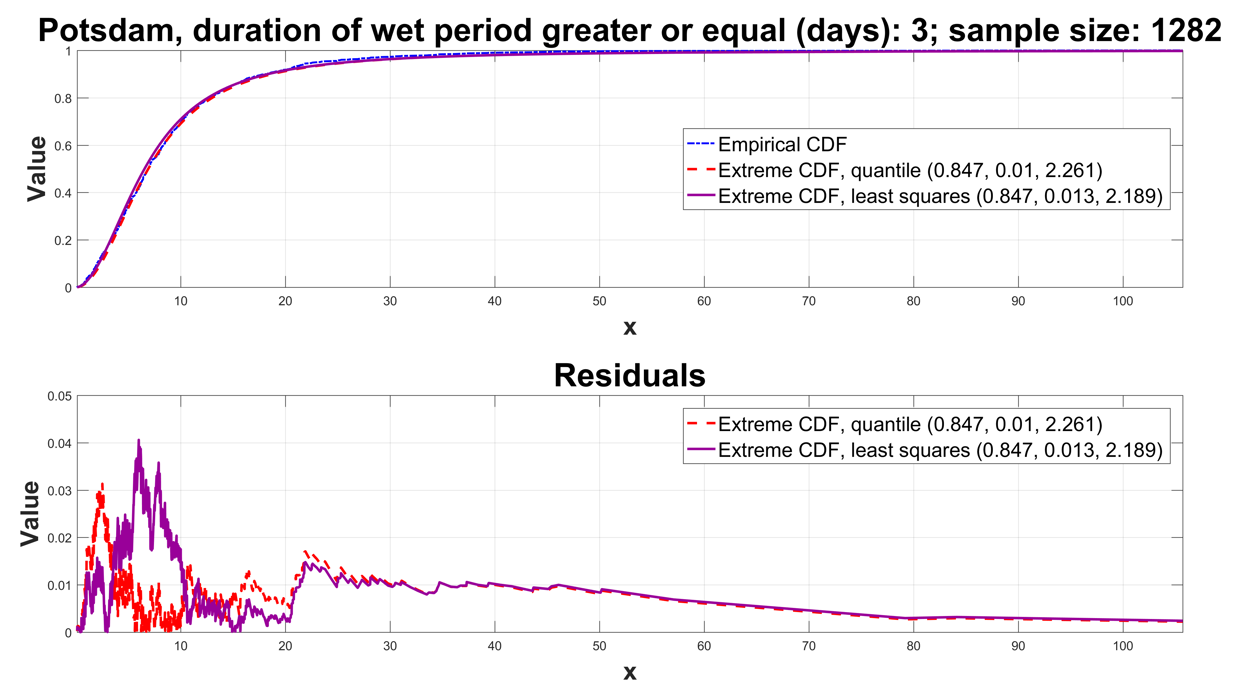}} c) \\
\end{minipage}
\hfill
\begin{minipage}[h]{0.44\linewidth}
\center{\includegraphics[width=1\linewidth, height=0.7\linewidth]{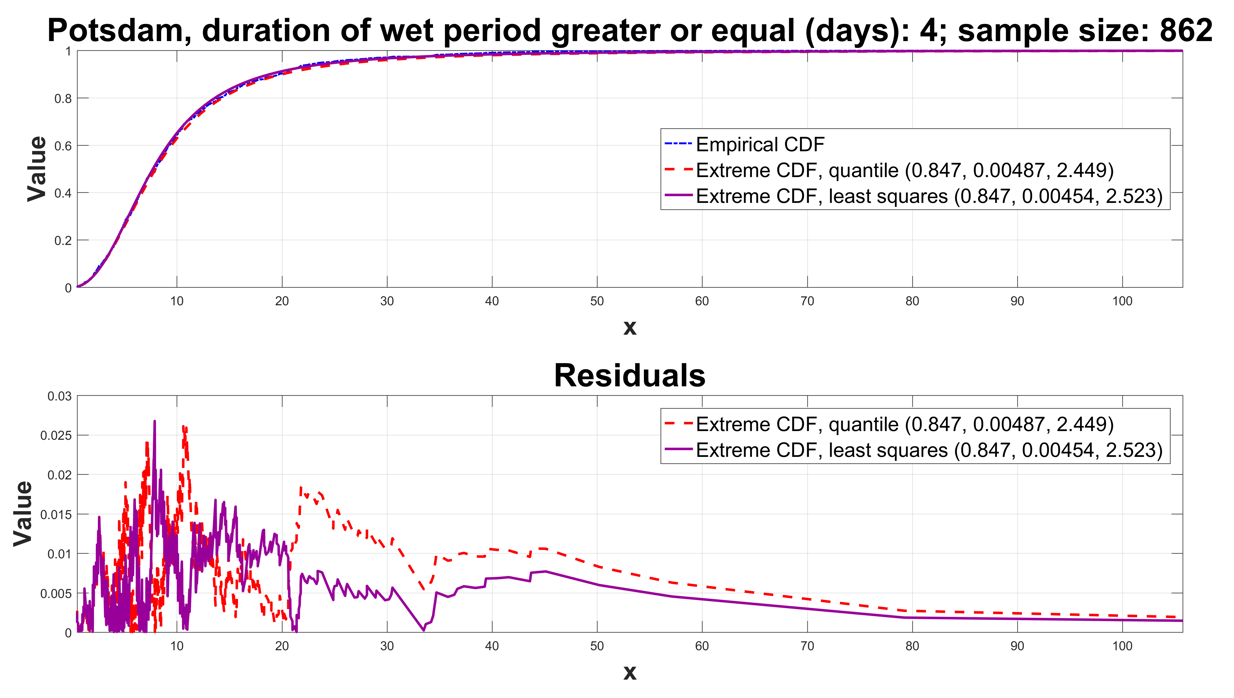}} d) \\
\end{minipage}
\hfill
\begin{minipage}[h]{0.44\linewidth}
\center{\includegraphics[width=1\linewidth, height=0.7\linewidth]{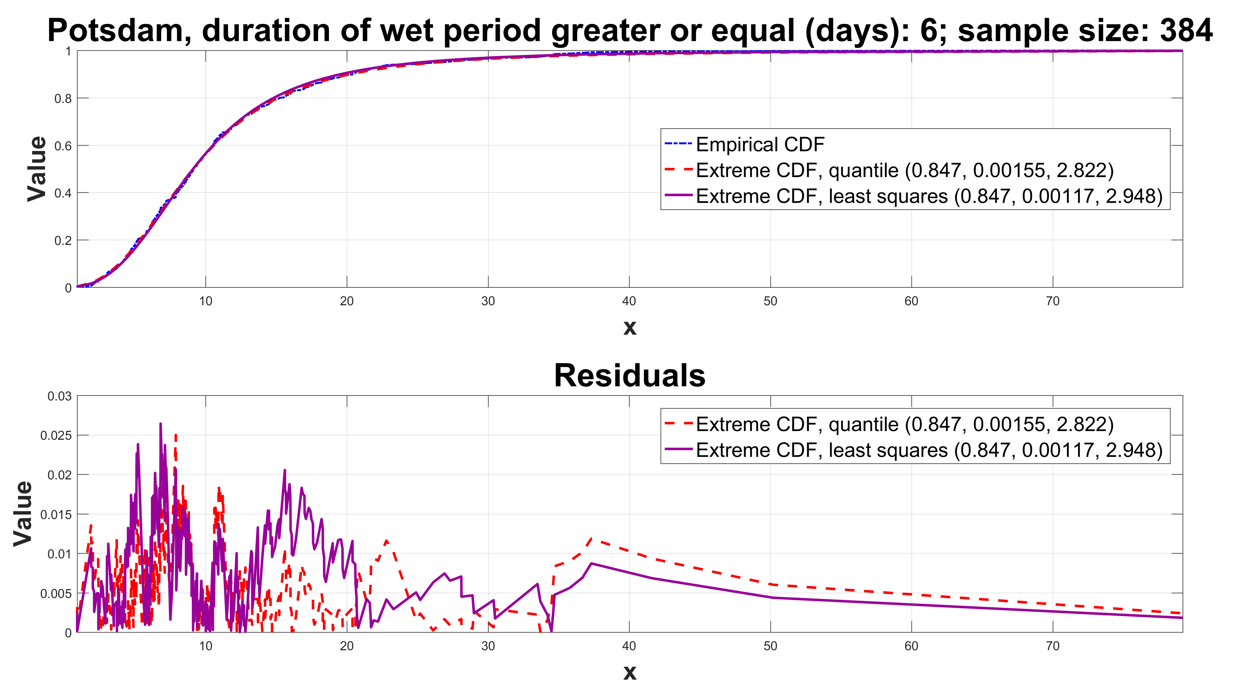}} e) \\
\end{minipage}
\hfill
\begin{minipage}[h]{0.44\linewidth}
\center{\includegraphics[width=1\linewidth, height=0.7\linewidth]{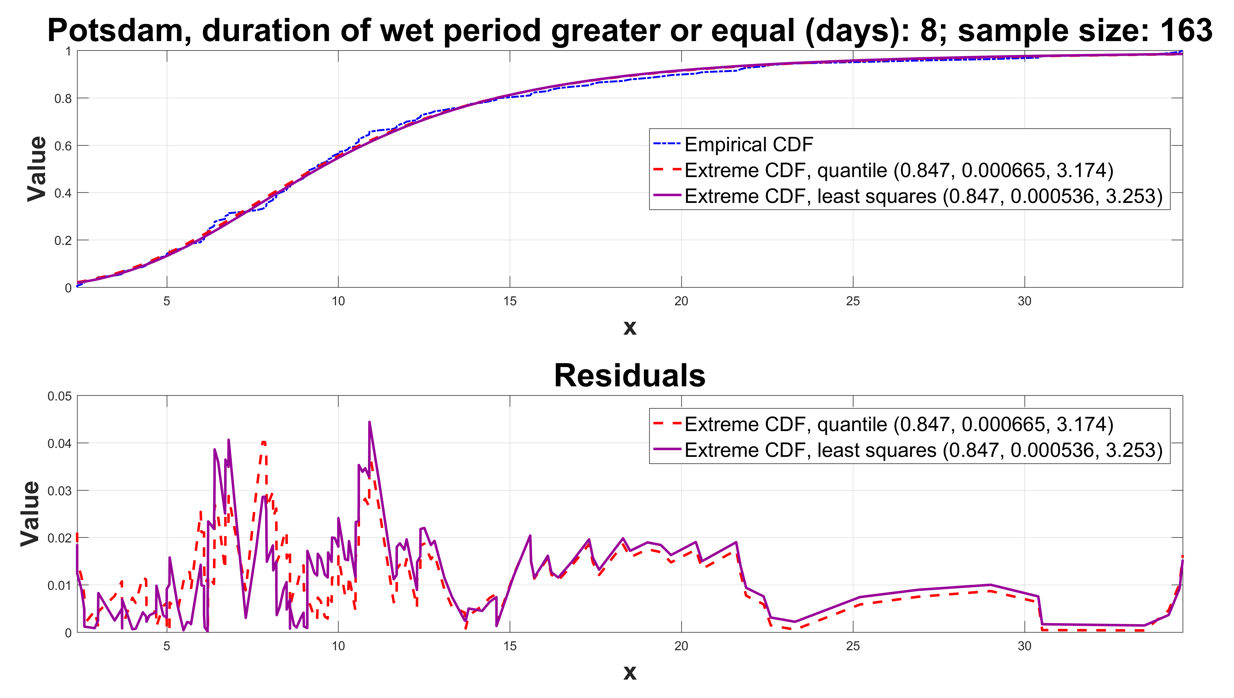}} f)\\
\end{minipage}
\vfill
\begin{minipage}[h]{0.44\linewidth}
\center{\includegraphics[width=1\linewidth, height=0.7\linewidth]{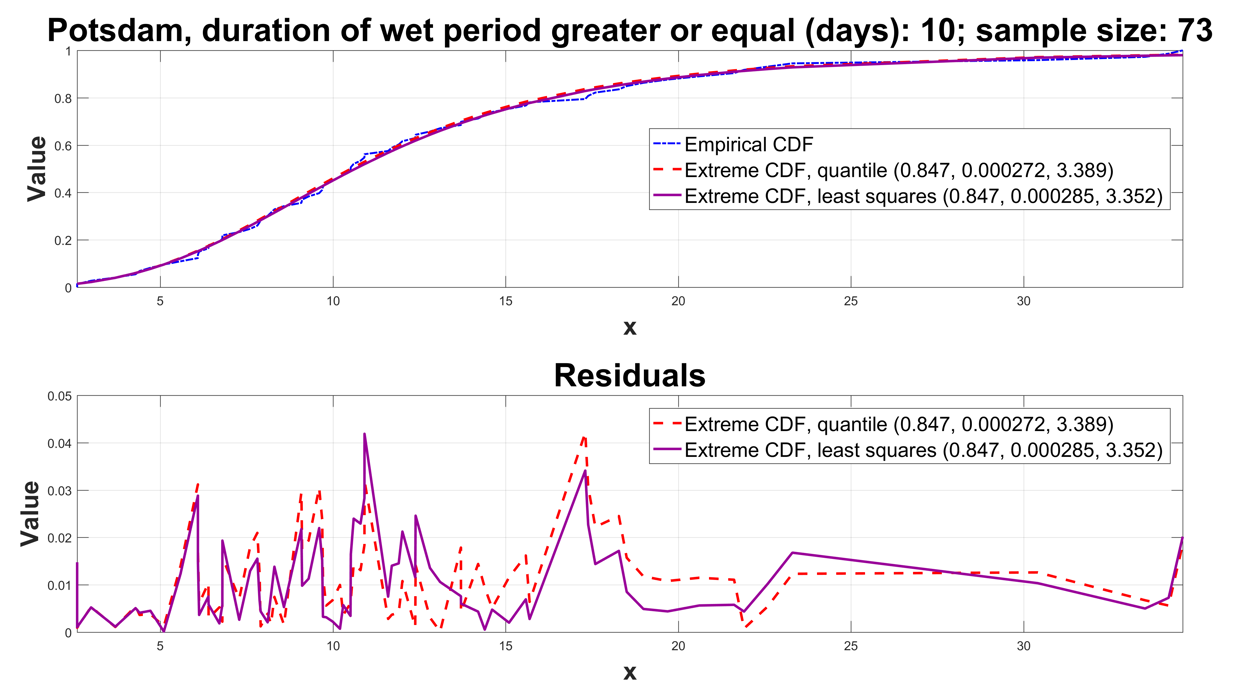}} g) \\
\end{minipage}
\hfill
\begin{minipage}[h]{0.44\linewidth}
\center{\includegraphics[width=1\linewidth, height=0.7\linewidth]{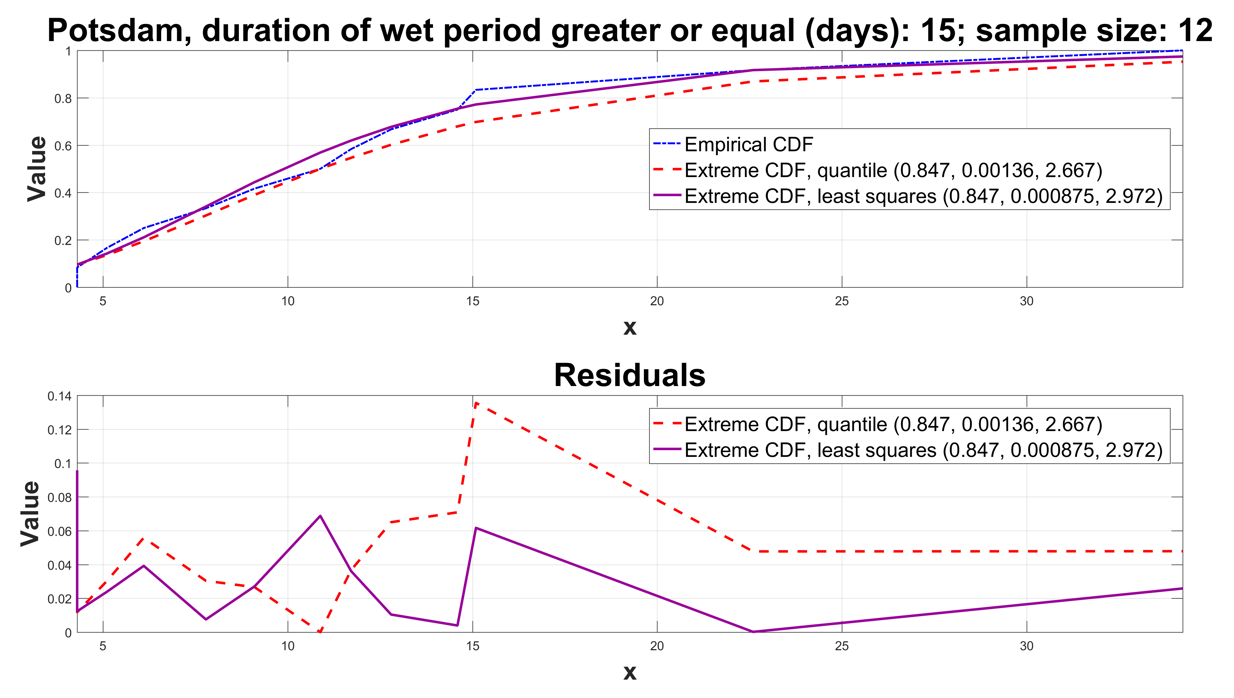}} h) \\
\end{minipage}
\caption{Empirical and estimated model extreme precipitation
distributions (Potsdam). Duration of wet period is no less than: a)
1; b) 2; c) 3; d) 4; e) 6; f) 8; g)~10; h)~15 days.}
\label{FigPotsdam}
\end{figure}

\begin{figure}[ptb]
\begin{minipage}[h]{0.44\linewidth}
\center{\includegraphics[width=1\linewidth, height=0.7\linewidth]{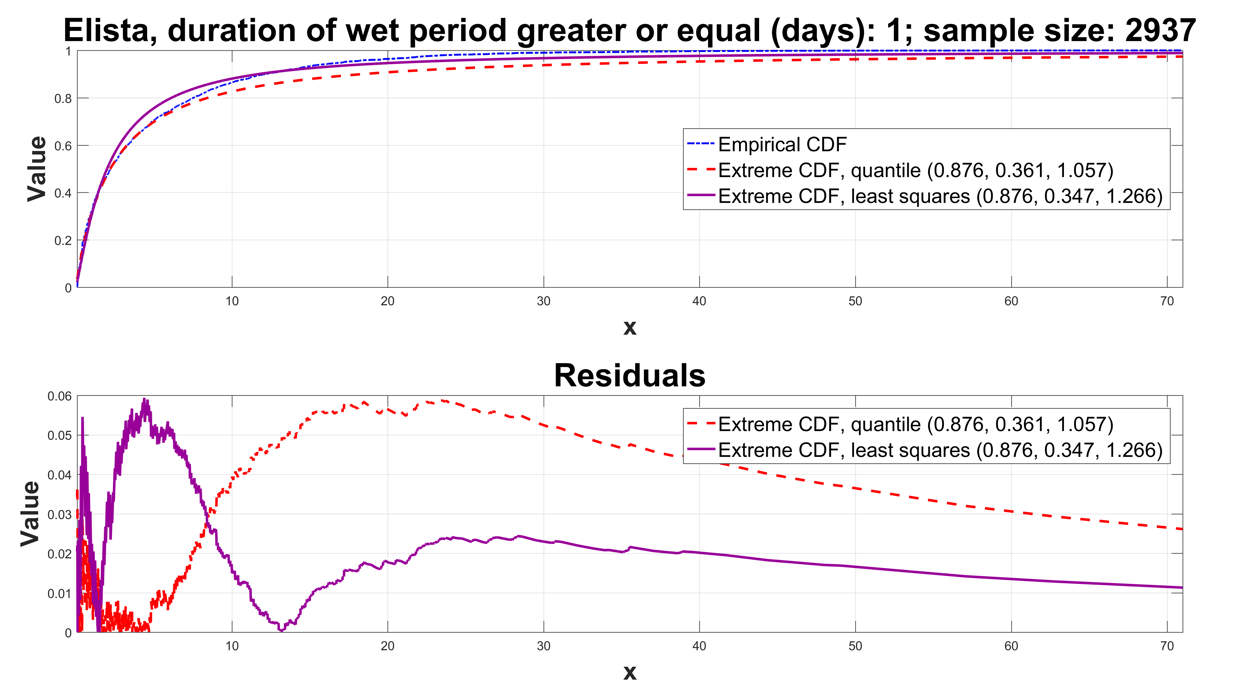}} a) \\
\end{minipage}
\hfill
\begin{minipage}[h]{0.44\linewidth}
\center{\includegraphics[width=1\linewidth,
height=0.7\linewidth]{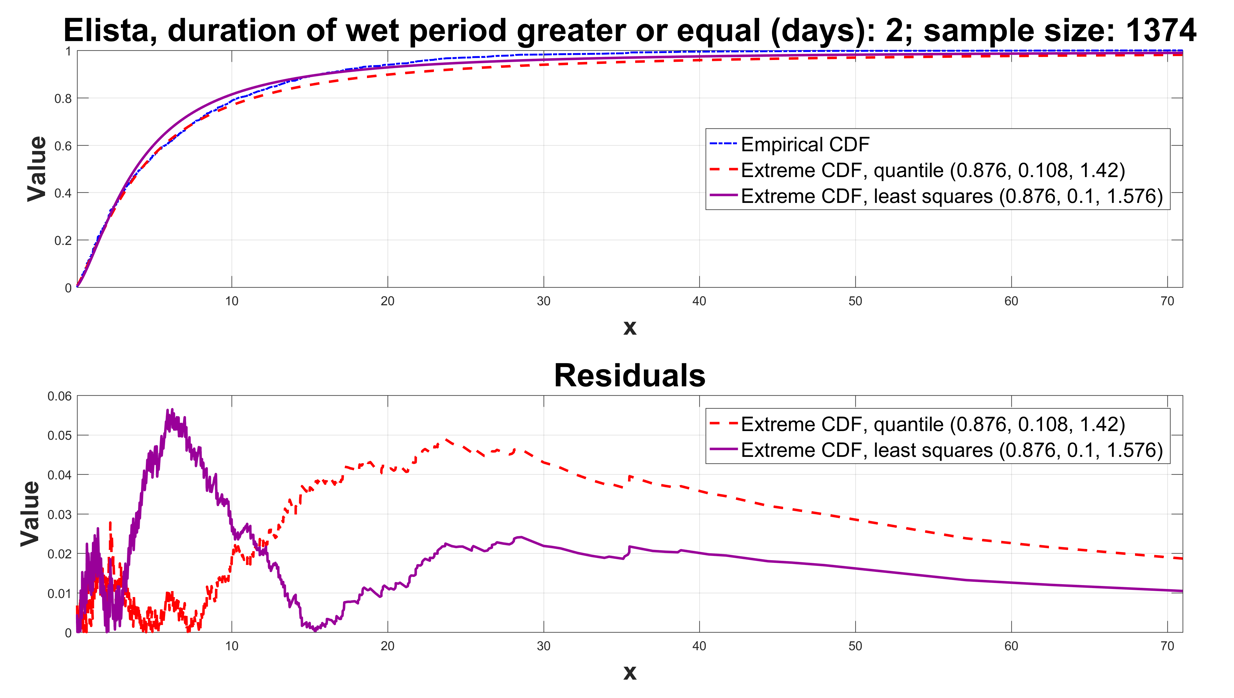}} b)\\
\end{minipage}
\vfill
\begin{minipage}[h]{0.44\linewidth}
\center{\includegraphics[width=1\linewidth, height=0.7\linewidth]{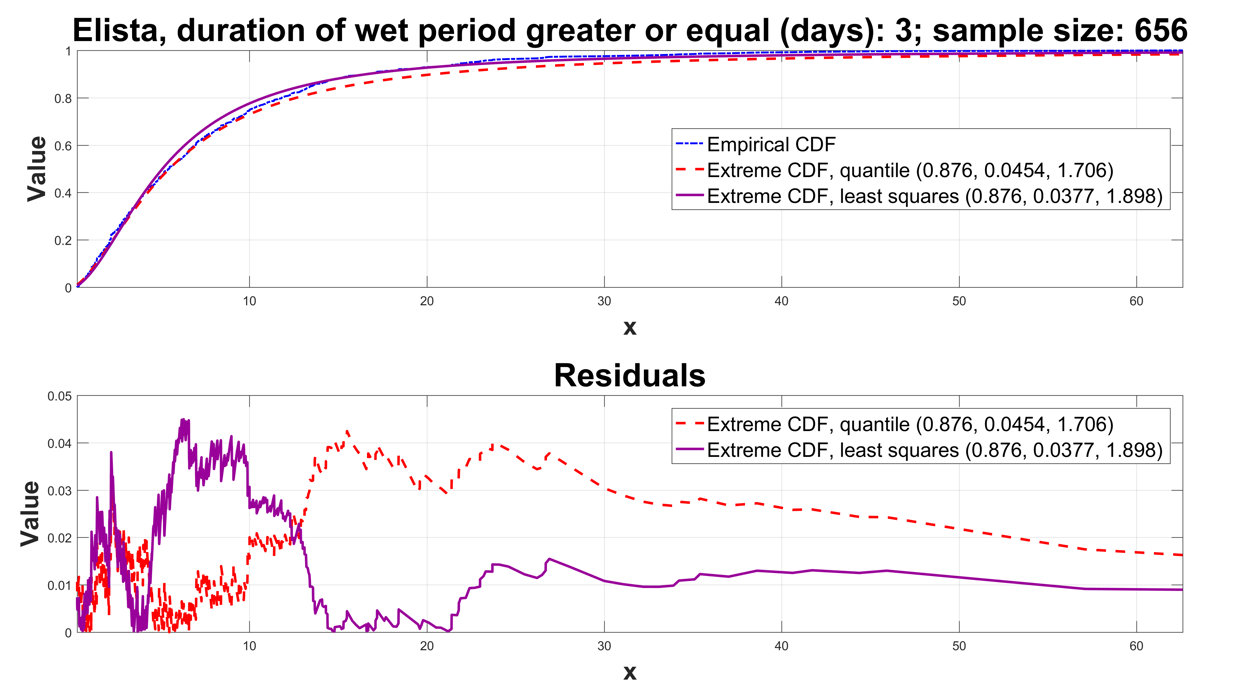}} c) \\
\end{minipage}
\hfill
\begin{minipage}[h]{0.44\linewidth}
\center{\includegraphics[width=1\linewidth, height=0.7\linewidth]{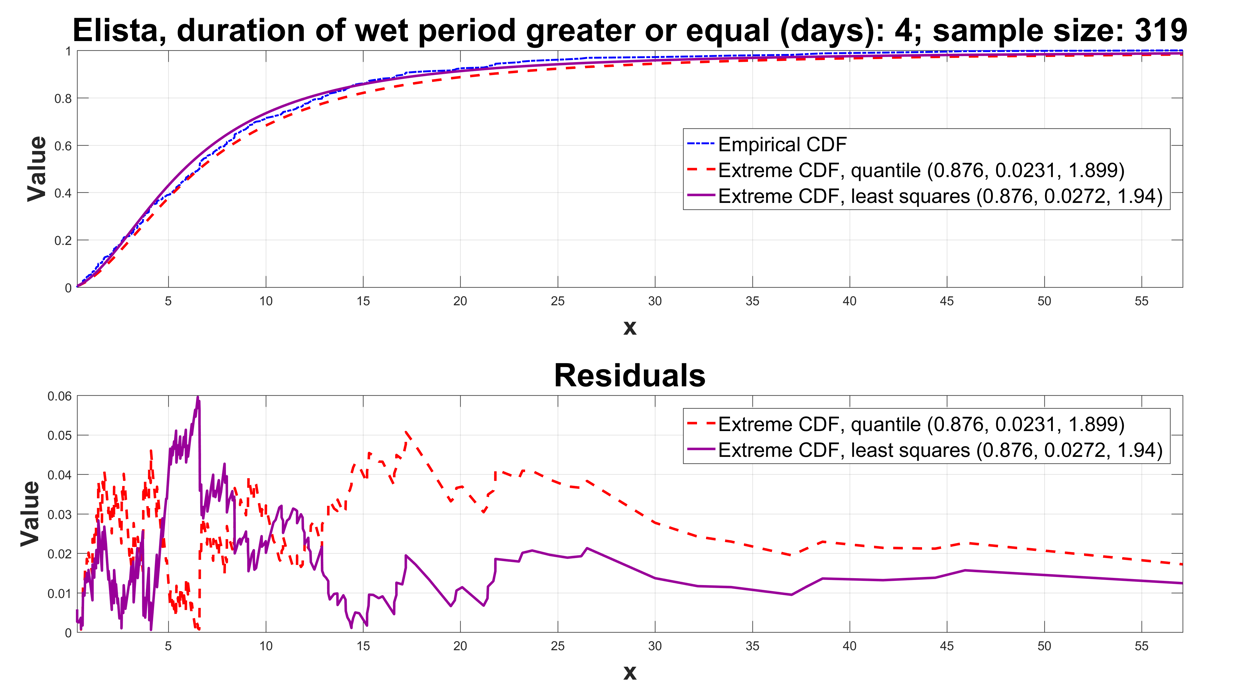}} d) \\
\end{minipage}
\hfill
\begin{minipage}[h]{0.44\linewidth}
\center{\includegraphics[width=1\linewidth, height=0.7\linewidth]{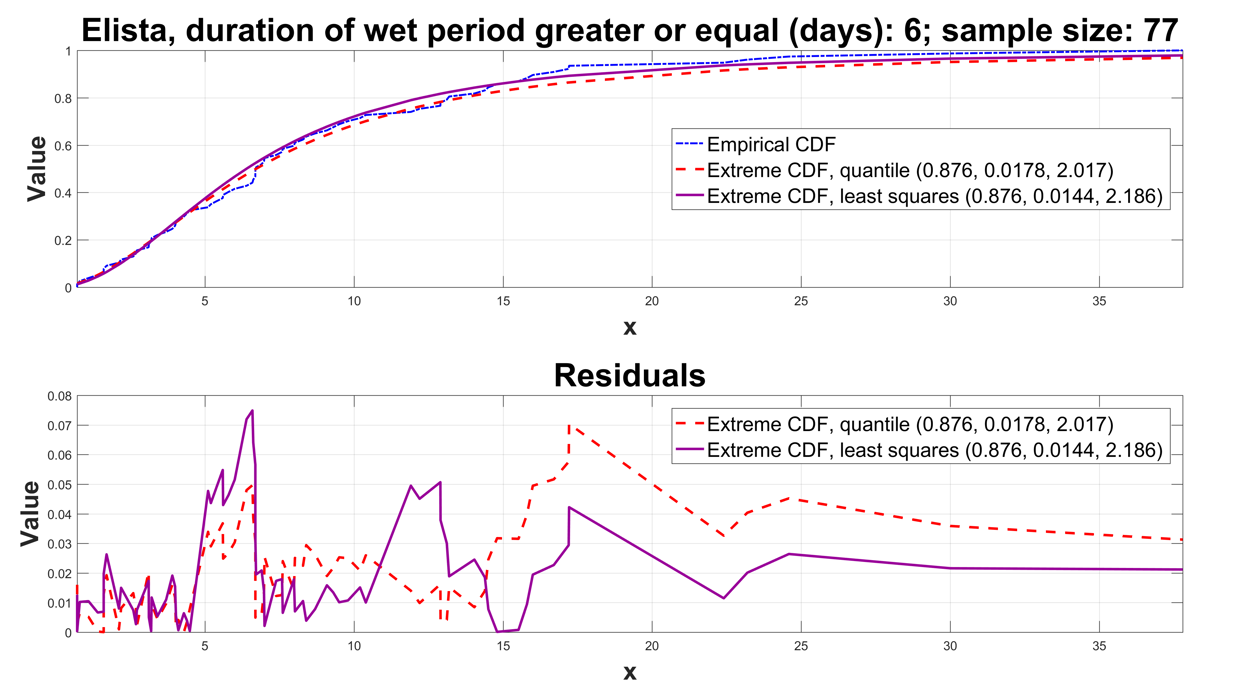}} e) \\
\end{minipage}
\hfill
\begin{minipage}[h]{0.44\linewidth}
\center{\includegraphics[width=1\linewidth, height=0.7\linewidth]{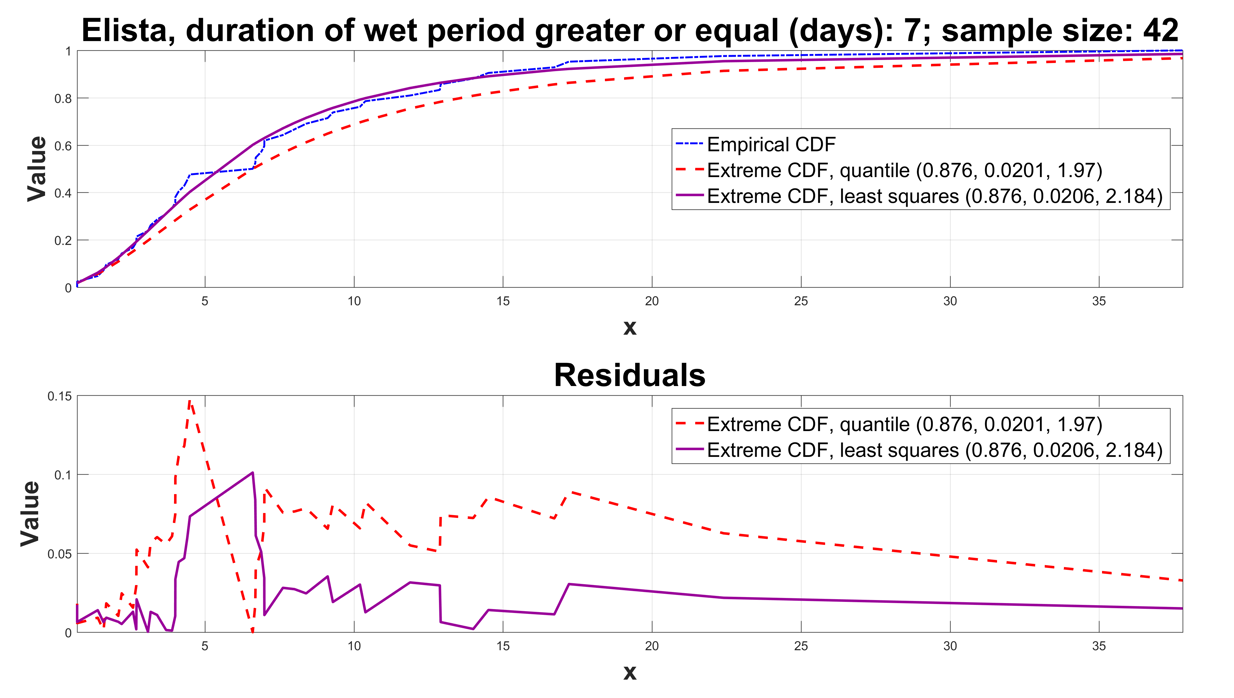}} f)\\
\end{minipage}
\vfill
\begin{minipage}[h]{0.44\linewidth}
\center{\includegraphics[width=1\linewidth, height=0.7\linewidth]{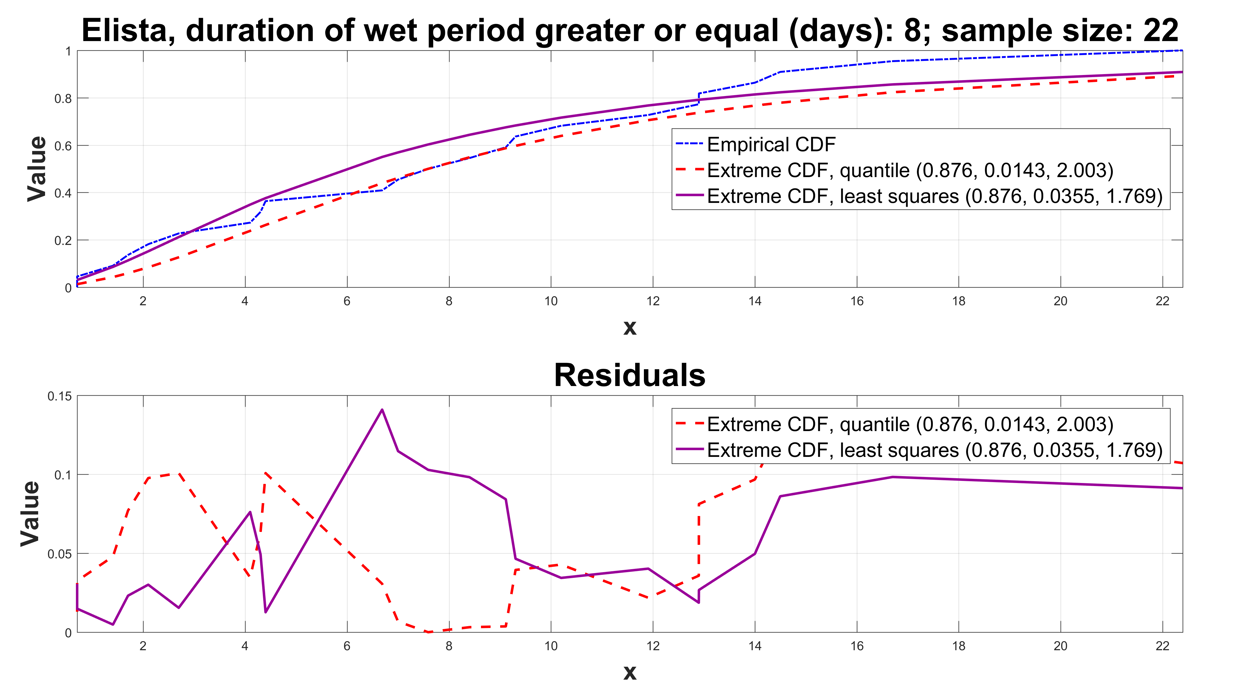}} g) \\
\end{minipage}
\hfill
\begin{minipage}[h]{0.44\linewidth}
\center{\includegraphics[width=1\linewidth, height=0.7\linewidth]{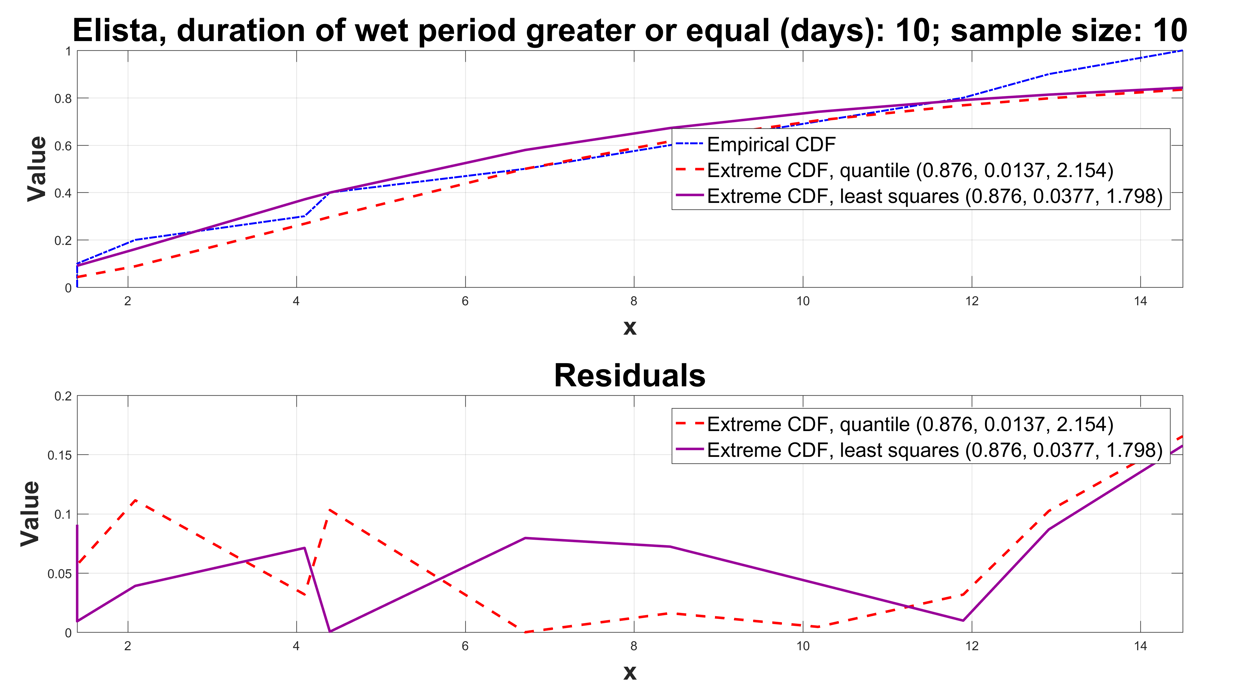}} h) \\
\end{minipage}
\caption{Empirical and estimated model extreme precipitation
distributions (Elista). Duration of wet period is no less than: a)
1; b) 2; c) 3; d) 4; e) 6; f) 7; g)~8; h) 10 days.}
\label{FigElista}
\end{figure}

The parameter $r$ is assumed known and its value coincides with that
of the negative binomial distribution (see the Introduction):
$r=0.847$ for Potsdam and $r=0.876$ for Elista.
Figures~\ref{FigPotsdam} and~\ref{FigElista} illustrate the
approximation of the empirical d.f. by the model d.f.
$F(x;r,\gamma,\lambda)$ with $\gamma$ and $\lambda$ estimated by the
`rough' formulas~\eqref{gamma} and~\eqref{lambda} as well as by the
least squares formulas~\eqref{gammaLS} and~\eqref{lambdaLS}. To
illustrate the asymptotic character of the approximation
$F(x;r,\gamma,\lambda)$ we consider a sort of censoring in which the
censoring threshold is the minimum length of the wet periods which
varies from 1 day (full sample) to 15 days.

For each censoring threshold $h=\min_im_i$ the sample is formed
according to the rule~\eqref{VarSample}.

For each value of the threshold on the upper graph there are
\begin{itemize}
\item the empirical d.f. (dot-dash line);
\item the d.f. $F(x;,r,\gamma,\lambda)$ with $\gamma$ and $\lambda$
estimated by the `rough' formulas~\eqref{gamma} and~\eqref{lambda}
(dash line);
\item the d.f. $F(x;,r,\gamma,\lambda)$ with $\gamma$ and $\lambda$
estimated by the least squares formulas~\eqref{gammaLS}
and~\eqref{lambdaLS} (continuous line).
\end{itemize}

On the lower graph there is the discrepancy (the uniform distance)
between the empirical d.f. and the fitted model d.f. The types of
lines correspond to those on the upper graph.

First of all, from Figures~\ref{FigPotsdam} and~\ref{FigElista} it
is seen that the asymptotic model $F(x;r,\gamma,\lambda)$ provides
very good approximation to the real statistical regularities in the
behavior of extremal daily precipitation within a wet period. As
this is so, the least squares formulas~\eqref{gammaLS}
and~\eqref{lambdaLS} yield more accurate estimates for the
parameters of the model d.f.

It should be also noted that these figures illustrate the dependence
of the accuracy of the approximation on the censoring threshold and
the censored sample size. The approximation is satisfactory even if
the censoring threshold $h$ is greater or equal to three and the
censored sample size is grater than 150. As this is so, the
influence of the threshold $h$ on the accuracy is more noticeable
than that of the sample size.

\section*{Acknowledgements}
The research was partly supported by the RAS Presidium Program No. I.33P (project 0063-2016-0015) and by
the Russian Foundation for Basic Research (projects 15-07-04040 and 17-07-00851).


\begin{thebibliography}{100}

\bibitem{Zolina2013} \emph{Zolina~O., Simmer~C., Belyaev~K., Gulev~S., Koltermann~P.}
Changes in the duration of European wet and dry spells during the last 60 years //
Journal of Climate, 2013. Vol. 26. P. 2022--2047.

\bibitem{Gulev} \emph{Korolev~V.\,Yu., Gorshenin~A.\,K., Gulev~S.\,K., Belyaev~K.\,P.,
Grusho~A.\,A.} Statistical Analysis of Precipitation Events~// AIP Conference Proceedings, 2017. Also available on arXiv:1705.11055 [math.PR]

\bibitem{Gorshenin2017} \emph{Gorhenin A.\,K.} On some mathematical and program methods for the construction of structure models of information flows~// Informatics
and Its Applications, 2017. Vol.~11. No.~1. P.~58--68.

\bibitem{Kingman1993} \emph{Kingman~J.\,F.\,C.} Poisson processes. -- Oxford: Clarendon Press,
1993.

\bibitem{GnedenkoKorolev1996} \emph{Gnedenko~B.\,V., Korolev~V.\,Yu.} Random Summation:
Limit Theorems and Applications. -- Boca Raton: CRC Press, 1996.

\bibitem{KorolevBeningShorgin2011} \emph{Korolev~V.\,Yu., Bening~V.\,E., Shorgin S.\,Ya.}
Mathematical Foundations of Risk Theory. 2nd Edition. -- Moscow:
FIZMATLIT, 2011.

\bibitem{Korolev2017} \emph{Korolev~V.\,Yu.} Analogs of the Gleser theorem for negative binomial and generalized gamma distributions and some their applications //
Informatics and Its Applications, 2017. Vol. 11. No. 3. 

\bibitem{Korolev2016TVP} \emph{Korolev~V.\,Yu.} Limit distributions for doubly stochastically
rarefied renewal processes and their properties // Theory Probab.
Appl., 2016. Vol.~61. No.~4. P.~753--773.

\bibitem{KorolevPoisson}
\emph{Korolev~V.\,Yu., Korchagin A.\,Yu., Zeifman A.\,I.} Poisson
theorem for the scheme of Bernoulli trials with random probability
of success and a discrete analog of the Weibull distribution //
Informatics and Its Applications, 2016. Vol. 10. No. 4. P. 11--20.

\bibitem{Korolev2016} \emph{Korolev~V.\,Yu., Korchagin~A.\,Yu., Zeifman~A.\,I.} On doubly
stochastic rarefaction of renewal processes // AIP Conference Proceedings, 2017. 

\bibitem{Gleser1989} \emph{Gleser~L.\,J.} The gamma distribution as a mixture of
exponential distributions // American Statistician, 1989. Vol. 43. P. 115--117.

\bibitem{Stacy1962} \emph{Stacy E. W.} A generalization of the gamma
distribution // Annals of Mathematical Statistics, 1962. Vol.~33.
P.~1187--1192.

\bibitem{KorolevZaks2013}
\emph{Zaks L.\,M., Korolev V.\,Yu.} Generalized variance gamma
distributions as limit laws for random sums // Informatics and Its
Applications, 2013. Vol. 7. No. 1. P. 105--115.

\bibitem{Zolotarev1983} \emph{Zolotarev V.\,M.} One-Dimensional Stable
Distributions. -- Providence, R.I.: American Mathematical Society,
1986.

\bibitem{Bolshev} \emph{Johnson, N.\,L., Kotz S., Balakrishnan N.} Continuous
Univariate Distributions, Vol. 2 (2nd Edition). -- New York: Wiley,
1995.

\bibitem{KorolevWeibull2016} \emph{Korolev V.\,Yu.} Product representations for random variables with the
Weibull distributions and their applications // Journal of
Mathematical Sciences, 2016. Vol. 218. No. 3. P. 298--313.

\bibitem{KotzOstrovskii1996} \emph{Kotz S., Ostrovskii I.\,V.} A mixture representation of the
Linnik distribution // Statistics and Probability Letters, 1996.
Vol. 26. P. 61--64.

\bibitem{KorolevZeifman2016a} \emph{ Korolev V.\,Yu., Zeifman A.\,I.} A note on mixture
representations for the Linnik and Mittag-Leffler distributions and
their applications // Journal of Mathematical Sciences, 2017. Vol.
218. No. 3. P. 314--327.

\bibitem{KorolevZeifman2016b} \emph{ Korolev V.\,Yu., Zeifman A.\,I.} Convergence of
statistics constructed from samples with random sizes to the Linnik
and Mittag-Leffler distributions and their generalizations //
Journal of Korean Statistical Society. Available online 25 July
2016. Also available on arXiv:1602.02480v1 [math.PR].

\bibitem{GreenwoodYule1920} \emph{Greenwood~M., Yule~G.\,U.} An inquiry into
the nature of frequency-distributions of multiple happenings, etc.
// J. Roy. Statist. Soc., 1920. Vol. 83. P. 255--279.

\bibitem{Korolev1998} \emph{Korolev V.\,Yu.}  On convergence of distributions of
compound Cox processes to stable laws // Theory of Probability and
its Applications, 1999. Vol. 43. Vol. 4. P. 644--650.

\bibitem{KorolevSokolov2008} \emph{Korolev V.\,Yu., Sokolov I.\,A.}
Mathematical Models of Inhomogeneous Flows of Extremal Events. --
Moscow: Torus-Press, 2008.

\bibitem{Goldie1967} \emph{Goldie C. M.} A class of infinitely divisible distributions //
Math. Proc. Cambridge Philos. Soc., 1967. Vol. 63. P. 1141-1143.

\end{thebibliography}
\end{document}